\documentclass[11pt]{article}
\usepackage{amssymb,latexsym,amsmath,amsbsy,amsthm,amsxtra,amsgen}
\newfont{\msbm}{msbm10 at 11pt}

\def\be{\begin{equation}}
\def\ee{\end{equation}}
\def\ba{\begin{align}}
\def\ea{\end{align}}

\newtheorem{Theo}{Theorem}
\newtheorem{Lemma}[Theo]{Lemma}
\newtheorem{Cor}[Theo]{Corollary}
\newtheorem{Prop}[Theo]{Proposition}

\begin{document}
\title{The waiting time for $m$ mutations}
\author{by Jason Schweinsberg\thanks{Supported in part by NSF Grant DMS-0504882}
\\ University of California, San Diego}
\maketitle



\begin{abstract}
We consider a model of a population of fixed size $N$ in which each individual gets replaced at rate one and each individual experiences a mutation at rate $\mu$.  We calculate the asymptotic distribution of the time that it takes before there is an individual in the population with $m$ mutations.  Several different behaviors are possible, depending on how $\mu$ changes with $N$.  These results have applications to the problem of determining the waiting time for regulatory sequences to appear and to models of cancer development.
\end{abstract}

\section{Introduction}
It is widely accepted that many types of cancer arise as a result of not one but several mutations.  For example, Moolgavkar and Luebeck \cite{ml92} write that ``the concept of multistage carcinogenesis is one of the central dogmas of cancer research", while Beerenwinkel et. al. \cite{beer07} write that ``the current view of cancer is that tumorigenesis is due to the accumulation of mutations in oncogenes, tumor suppressor genes, and genetic instability genes."  The idea that several mutations are required for cancer goes back at least to 1951, when Muller \cite{muller} wrote, ``There are, however, reasons for inferring that many or most cancerous growths would require a series of mutations in order for cells to depart sufficiently from the normal."  Three years later, Armitage and Doll \cite{armdoll} proposed a simple mathematical multi-stage model of cancer.  Motivated by the goal of explaining the power law relationship between age and incidence of cancer that had been observed by Fisher and Holloman \cite{fh51} and Nordling \cite{n53}, they formulated a model in which a cell that has already experienced $k-1$ mutations experiences a $k$th mutation at rate $u_k$.  They showed that asymptotically as $t \rightarrow 0$, the probability that the $m$th mutation occurs in the time interval $[t, t + dt]$ is given by
\begin{equation} \label{adrate}
r(t) \: dt = \frac{u_1 u_2 \dots u_m t^{m-1}}{(m-1)!} \: dt.
\end{equation}
They fit their model to data from 17 different types of cancer, and found that for many types of cancer the incidence rate $r(t)$ increases like the fifth or sixth power of age, suggesting that perhaps 6 or 7 mutations are involved in cancer progression.  Because of concerns that having 6 or 7 stages may not be biologially plausible, Armitage and Doll \cite{armdoll57} later proposed a two-stage model as an alternative.  A more general two-stage model was proposed by Moolgavkar and Knudson \cite{mknud81}, who demonstrated that two-stage models are flexible enough to fit a wide range of data if one allows for the possibilities that the number of healthy cells with no mutations may change over time, and that cells with one mutation may divide rapidly, causing the second mutation, and therefore the onset of cancer, to happen more quickly than it otherwise would.

Since the seminal papers of Armitage and Doll, multi-stage models have been applied to a number of different types of cancer.  Knudson \cite{knud71, hknud78} discovered that retinoblastoma is a result of getting two mutations.  Multi-stage models of colon cancer have been studied extensively.  Moolgavkar and Luebeck \cite{ml92} argued that a three-stage model fit the available data slightly better than a two-stage model.  Later in \cite{lm02}, they found a good fit to a four-stage model.  Calabrese et. al. \cite{cal05} worked with data from 1022 cancers from 9 hospitals in Finland and estimated that between 4 and 9 mutations are required for cancer, with fewer mutations being required for hereditary cancers than for sporadic (nonhereditary) cancers.  A recent study \cite{sj06} of over 13,000 genes from breast and colon cancers suggests that as many as 14 mutations may be involved in colon cancer and as many as 20 may be involved in breast cancer.  Multi-stage models have also been fit to data on lung cancer \cite{fn89} and T-cell leukemia \cite{k90}.  See \cite{knud01} for a recent survey of applications of multi-stage cancer models.

In this paper, we formulate a simple mathematical model and calculate the asymptotic distribution of the time that it takes for cancer to develop.  Our model is as follows.  Consider a population of fixed size $N$.  We think of the individuals in the population as representing $N$ cells, which could develop cancer.  We assume that the population evolves according to the Moran model \cite{moran}.  That is, each individual independently lives for an exponentially distributed amount of time with mean one, and then is replaced by a new individual whose parent is chosen at random from the $N$ individuals in the population (including the one being replaced).  These births and deaths represent cell division and cell death.  We also assume that each individual independently experiences mutations at times of a rate $\mu$ Poisson process, and each new individual born has the same number of mutations as its parent.  We refer to an individual that has $j$ mutations as a type $j$ individual, and a mutation that takes an individual's number of mutations from $j-1$ to $j$ as a type $j$ mutation.  Let $X_j(t)$ be the number of type $j$ individuals at time $t$.  For each positive integer $m$, let $\tau_m = \inf\{t: X_m(t) > 0\}$ be the first time at which there is an individual in the population with $m$ mutations.  We view $\tau_m$ as representing the time that it takes for cancer to develop.  Clearly $\tau_1$ has the exponential distribution with rate $N \mu$ because the $N$ individuals are each experiencing mutations at rate $\mu$.  Our goal in this paper is to compute the asymptotic distribution of $\tau_m$ for $m \geq 2$.

When a new mutation occurs, eventually either all individuals having the mutation die, causing the mutation to disappear from the population, or the mutation spreads to all individuals in the population, an event which we call fixation.  Because a mutation initially appears on only one individual and is assumed to offer no selective advantage or disadvantage, each mutation fixates with probability $1/N$.  
Once one mutation fixates, the problem reduces to waiting for $m-1$ additional mutations.  However, it is possible for one individual to accumulate $m$ mutations before any mutation fixates in the population, an event which is sometimes called stochastic tunneling (see \cite{imn04}).  It is also possible for there to be $j$ fixations, and then for one individual to get $m-j$ mutations that do not fixate.  Because there are different ways to get $m$ mutations, the limiting behavior is surprisingly complex, as the form of the limiting distribution of $\tau_m$ depends on how $\mu$ varies as a function of $N$.

There is another source of biological motivation for this model coming from the evolution of regulatory sequences.  Regulatory sequences are short DNA sequences that control how genes are expressed.  Getting a particular regulatory sequence would require several mutations, so to understand the role that regulatory sequences play in evolution, one needs to understand how long it takes before these mutations occur.  See Durrett and Schmidt \cite{regseq1, regseq2} for work in this direction.

In addition to this motivation from biology, there is mathematical motivation for studying this model as well.  The model is simple and natural and, as will be seen from the results, gives rise to different asymptotic behavior depending on how $\mu$ scales as a function of $N$.  In particular, the usual diffusion scaling from population genetics in which $N \mu$ tends to a constant is just one of several regimes.

This paper can be viewed as a sequel to \cite{DSS}, in which the authors considered a more general model in which an individual with $k-1$ mutations experiences a $k$th mutation at rate $u_k$.  The model considered here is the special case in which $u_k = \mu$ for all $k$, so we are assuming that all mutation rates are the same.  However, whereas in \cite{DSS} results were obtained only for specific ranges of the mutation rates $u_k$, here we are able to obtain all possible limiting behaviors for the case in which the mutation rates are the same.  We also emphasize that although our model accounts for cell division and cell death, we assume that the rates of cell division and cell death are the same, unlike many models in the biology literature which specify that individuals with between 1 and $m-1$ mutations have a selective advantage, allowing their numbers to increase rapidly (see, for example, \cite{armdoll57, mknud81, ml90, ml92, beer07}).  As we explain below, several special cases of our results have previously appeared in the biology literature, especially for the two-stage models when $m = 2$.  However, here we are able to give complete asymptotic results for all $m$, as well as to provide rigorous proofs of the results.  We state our main results in section 2.  Proofs are given in sections 3, 4, and 5.

\section{Main results}

In this section, we state our results on the limiting behavior of the waiting time for an individual to acquire $m$ mutations, and we explain the heuristics behind the results.  Many of the heuristics are based on approximation by branching processes.  In the Moran model, if $k$ individuals have a mutation, then the number of individuals with the mutation is decreasing by one at rate $k(N-k)/N$ (because the $k$ individuals with the mutation are dying at rate $k$, and the probability that the replacement individual does not have a mutation is $(N-k)/N$) and is increasing by one at rate $k(N-k)/N$ (because the $N-k$ individuals without a mutation are dying at rate one, and the replacement individual has a mutation with probability $k/N$).  Therefore, when $k$ is much smaller than $N$, the number of individuals with a given mutation behaves approximately like a continuous-time branching process in which each individual gives birth and dies at rate one.

To keep track of further mutations, it is natural to consider a continuous-time multitype branching process in which initially there is a single type 1 individual, each individual gives birth and dies at rate 1, and a type $j$ individual mutates to type $j+1$ at rate $\mu$.  If $p_j$ denotes the probability that there is eventually a type $j$ individual in the population, then
\begin{equation}\label{pjeq}
p_j = \frac{1}{2 + \mu}(2p_j - p_j^2) + \frac{\mu}{2+\mu} p_{j-1}.
\end{equation}
To see this result, condition on the first event.  With probability $1/(2 + \mu)$, the first event is a death, and there is no chance of getting a type $j$ individual.  With probability $1/(2 + \mu)$, the first event is a birth, in which case each individual has a type $j$ descendant with probability $p_j$ and therefore the probability that at least one has a type $j$ descendant is $2p_j - p_j^2$.  With probability $\mu/(2 + \mu)$, the first event is a mutation to type $2$, in which case the probability of a type $j$ descendant is $p_{j-1}$ because $j-1$ further mutations are needed.  Equation (\ref{pjeq}) can be rewritten as $p_j^2 + \mu p_j - \mu p_{j-1} = 0$, and the positive solution is $$p_j = \frac{-\mu + \sqrt{\mu^2 + 4 \mu p_{j-1}}}{2}.$$  When $\mu$ is small, the second term under the square root dominates the numerator, and we get $p_j \approx \sqrt{\mu p_{j-1}}$.  Since $p_1 = 1$, the approximation $p_j \approx \mu^{1 - 2^{-(j-1)}}$ follows by induction.

Because the Moran model can be approximated by a branching process when the number of mutant individuals is much smaller than $N$, this result suggests that under appropriate conditions, the probability that a type $1$ individual in the population has a type $m$ descendant should be approximately $\mu^{1 - 2^{-(m-1)}}$.  Proposition \ref{family} below, which is a special case of Proposition 4.1 in \cite{DSS}, establishes that this approximation is indeed valid.  Here and throughout the paper, the mutation rate $\mu$ depends on $N$ even though we do not record this dependence in the notation.  Also, if $f$ and $g$ are two functions of $N$, we write $f(N) \sim g(N)$ if $f(N)/g(N) \rightarrow 1$ as $N \rightarrow \infty$.  We also write $f(N) \ll g(N)$ if $f(N)/g(N) \rightarrow 0$ as $N \rightarrow \infty$ and $f(N) \gg g(N)$ if $f(N)/g(N) \rightarrow \infty$ as $N \rightarrow \infty$.

\begin{Prop}\label{family}
Consider a model which is identical to the model described in the introduction, except that initially there is one individual of type 1 and $N-1$ individuals of type 0, and no further type 1 mutations are possible.  Let $q_m$ be the probability that a type $m$ individual eventually is born.  Suppose that $N \mu^{1 - 2^{-(m-1)}} \rightarrow \infty$ as $N \rightarrow \infty$, and that there is a constant $a > 0$ such that $N^a \mu \rightarrow 0$.  Then $$q_m \sim \mu^{1 - 2^{-(m-1)}}.$$
\end{Prop}

Note that $q_m$ is the probability that a given type 1 individual eventually has a type $m$ descendant.  Because a number of our arguments involve considering each type 1 mutation and its descendants separately from other type 1 mutations, this result will be used repeatedly.

To understand the order of magnitude of $q_m$ another way, recall that the probability that the total progeny of a critical branching process exceeds $M$ is of order $M^{-1/2}$ (see, for example, \cite{harris}), so if there are $L$ independent branching processes, the most successful will have a total progeny of order $L^2$.  Furthermore, the sum of the total progenies of the $L$ processes will also be of order $L^2$.  Therefore, if there are $L$ type 1 mutations, the number of descendants they produce will be of order $L^2$.  Each type 1 descendant will experience a type 2 mutation before dying with probability approximately $\mu$, so this should lead to on the order of $L^2 \mu$ type 2 mutations.  It follows that the number of type 2 descendants should be on the order of $L^4 \mu^2$, and this will lead to on the order of $L^4 \mu^3$ type 3 mutations.  Repeating this reasoning, we see that the number of type $m$ mutations should be of order $L^{2^{m-1}} \mu^{2^{m-1}-1}$.  By setting this expression equal to one and solving for $L$, we see that it should take on the order of $\mu^{-(1-2^{-(m-1)})}$ type 1 mutations before one of these mutations gets a type $m$ descendant.  That is, the probability that a type 1 individual has a type $m$ descendant is of order $\mu^{1-2^{-(m-1)}}$. 

\subsection{Gamma limits when $N \mu \rightarrow 0$}

Because mutations occur at times of a Poisson process of rate $N \mu$, there will be approximately $N \mu T$ mutations by time $T$.  We have seen that after a mutation occurs, the number of individuals with the mutation behaves approximately like a critical branching process.  By a famous result of Kolmogorov \cite{kolm38}, the probability that a critical branching process survives for time $t$ is of order $1/t$.  This means that if we have $N \mu T$ independent critical branching processes, the most successful will survive for a time which is of order $N \mu T$.  Therefore, all mutations that appear before time $T$ should either die out or fixate after being in the population for a time of order $N \mu T$.  If $N \mu \ll 1$, then this time is much smaller than the time $T$ that we have to wait for the mutation.  Therefore, when $N \mu \ll 1$, we can consider each mutation separately and determine whether either it fixates or gives birth to a type $m$ descendant without fixating.  We can ignore the time that elapses between when the original mutation appears, and when either it fixates or the descendant with $m$ mutations is born.  The importance of the condition $N \mu \ll 1$ was previously noted, for example, in \cite{imn04} and \cite{ksn03}.

We have already seen that a mutation fixates with probability $1/N$ and gives birth to a type $j$ descendant with probability approximately $\mu^{1-2^{-(j-1)}}$.  Therefore, fixation of some mutation will happen first if $N \mu^{1-2^{-(j-1)}} \rightarrow 0$ as $N \rightarrow \infty$ or, equivalently, if $\mu \ll N^{-2^{j-1}/(2^{j-1}-1)}$.  This leads to the following result when $N \mu \ll 1$.  Note that when $m = 2$, the result in part 1 of the theorem matches (12.12) of \cite{nowak}, while the result in part 3 matches (12.14) of \cite{nowak}; see also section 3 of \cite{ksn03}, section 4 of \cite{imkn05}, and Theorem 1 of \cite{regseq2}.

\begin{Theo}\label{smallmu}
Let $Z_1, Z_2, \dots$ be independent random variables having the exponential distribution with rate $1$, and let $S_k = Z_1 + \dots + Z_k$, which has a gamma distribution with parameters $(k, 1)$.
\begin{enumerate}
\item If $\mu \ll N^{-2}$, then $\mu \tau_m \rightarrow_d S_{m-1}$.

\item If $N^{-2^{j-1}/(2^{j-1}-1)} \ll \mu \ll N^{-2^j/(2^j - 1)}$ for some $j = 2, \dots, m-1$, then $\mu \tau_m \rightarrow_d S_{m-j}$.

\item If $N^{-2^{m-1}/(2^{m-1}-1)} \ll \mu \ll N^{-1}$, then $N \mu^{2 - 2^{-(m-1)}} \tau_m \rightarrow_d Z_1$.
\end{enumerate}
\end{Theo}

To understand this result, note that in part 1 of the theorem, when $\mu \ll N^{-2}$, fixation occurs before any individual gets two mutations without a fixation.  Therefore, to get $m$ mutations, we have to wait for $m-1$ different mutations to fixate, and this is the sum of $m-1$ independent exponential waiting times.  The exponential random variables have rate parameter $\mu$, because there are mutations at rate $N \mu$ and each fixates with probability $1/N$, so mutations that fixate occur at rate $\mu$.
Once $m-1$ fixations have occurred, the $m$th mutation occurs quickly, at rate $N \mu$ rather than at rate $\mu$, so only the waiting times for the $m-1$ fixations contribute to the limiting distribution.  For part 2 of the theorem, when $N^{-2^{j-1}/(2^{j-1}-1)} \ll \mu \ll N^{-2^j/(2^j - 1)}$ for some $j = 2, \dots, m-1$, fixation occurs before an individual can accumulate $j+1$ mutations, but an individual can accumulate $j$ mutations before fixation.  Therefore, we wait for $m-j$ fixations, and then the remaining $j$ mutations happen without fixation.  Because the $j$ mutations without fixation happen on a faster time scale, the limit is a sum of $m-j$ exponential random variables.  In part 3, we get $m$ mutations before the first fixation, and there is an exponential waiting time until the first mutation that is successful enough to produce an offspring with  $m$ mutations.  Mutations happen at rate $N \mu$, and mutations are successful with probability approximately $\mu^{1-2^{-(m-1)}}$, which explains the time-scaling factor of $N \mu^{2 - 2^{-(m-1)}}$.

Part 3 of Theorem \ref{smallmu} is the special case of Theorem 2 of \cite{DSS} in which $u_j = \mu$ for all $j$.  Condition $(i)$ of that theorem becomes the condition $\mu \ll N^{-1}$, while condition $(iv)$ becomes the condition $N^{-2^{m-1}/(2^{m-1}-1)} \ll \mu$.  Parts 1 and 2 of Theorem \ref{smallmu} above are proved in section \ref{smallpf}.

\subsection{The borderline cases}

Theorem \ref{smallmu} does not cover the cases when $\mu$ is of the order $N^{-2^{j-1}/(2^{j-1} - 1)}$ for some $j$.  On this time scale, for the reasons discussed in the previous section, we can still neglect the time between when a mutation first appears in the population and when it either fixates or dies out because this time will be much shorter than the time we had to wait for the mutation to occur.  However, fixations happen on the same time scale as events in which an individual gets $j$ mutations without fixation.  Therefore,  to get to $m$ mutations, we start with $m - j$ fixations.  Then we can either have another fixation (followed by $j-1$ additional mutations, which happen on a faster time scale) or we can get $j$ mutations without any fixation.  The waiting time is the sum of $m-j$ independent exponential random variables with rate $\mu$ and another exponential random variable having the faster rate $\lambda_j \mu$.  The last exponential random variable comes from waiting for a mutation that either fixates or has a descendant with $j-1$ additional mutations but does not fixate.  This leads to the following result.

\begin{Theo}\label{border}
Suppose $\mu \sim A N^{-2^{j-1}/(2^{j-1} - 1)}$ for some $j = 2, \dots, m$ and some constant $A > 0$.
Let $Z_1, Z_2, \dots$ be independent exponential random variables having the exponential distribution with rate $1$, and let $S_k = Z_1 + \dots + Z_k$.  Let $Y$ be independent of $Z_1, Z_2, \dots$, and assume that $Y$ has the exponential distribution with rate $\lambda_j$, where
\begin{equation} \label{lamj}
\lambda_j = \sum_{k=1}^{\infty} \frac{A^{2k(1 - 2^{-(j-1)})}}{(k-1)!(k-1)!} \bigg/ \sum_{k=1}^{\infty} \frac{A^{2k(1 - 2^{-(j-1)})}}{k!(k-1)!}.
\end{equation}
Then $\mu \tau_m \rightarrow_d S_{m-j} + Y$.
\end{Theo}

This result when $j = m$ is the special case of Theorem 3 of \cite{DSS} in which $u_j = \mu$ for all $j$.  As will be seen in section \ref{smallpf}, the result for $j \leq m-1$ follows easily from the result when $j = m$.

To explain where the formula for $\lambda_j$ comes from, we review here the outline of the proof of Theorem 3 in \cite{DSS}.  Assume that we already have $m-j$ fixations, and now we need to wait either for another fixation or for a mutation that will have a descendant with $j-1$ additional mutations.  We can not approximate the probability of the latter event by $\mu^{1-2^{-(j-1)}}$ in this case because to get $j-1$ further mutations, the number of individuals with the original mutation will need to be of order $N$, so the branching process approximation does not hold.  Instead, we consider a model in which there is one individual with a mutation at time zero, and $X(t)$ denotes the number of individuals with the mutation at time $t$.  At time $t$, the individuals with the mutation each experience further mutations at rate $\mu$, and these further mutations each have probability approximately $\mu^{1-2^{-(j-2)}}$ of having an offspring with $j$ total mutations.  Therefore, at time $t$, successful mutations are happening at rate $\gamma X(t)$, where $$\gamma \approx \mu \cdot \mu^{1-2^{-(j-2)}} = \mu^{2(1-2^{-(j-1)})}.$$  At time $t$, the jump rate of the process is $2 X(t)(N - X(t))/N$.  Therefore, by making a time-change, we can work instead with a continuous-time simple random walk $(Y(t), t \geq 0)$ which jumps at rate one, and the mutation rate at time $t$ becomes $$\gamma Y(t) \cdot \frac{N}{2 Y(t)(N - Y(t))} = \frac{\gamma}{2(1 - Y(t)/N)}.$$  Therefore, the probability that there is no fixation and no further successful mutation is approximately $$E \bigg[ \exp \bigg( - \int_0^T \frac{\gamma}{2(1 - Y(t)/N)} \: dt \bigg) {\bf 1}_{\{Y(T) = 0\}} \bigg],$$ where $T = \inf\{t: Y(t) \in \{0, N\} \}$.  Simple random walk converges to Brownian motion, so if instead of starting with just one mutant individual we assume that $Y(0) = \lfloor N x \rfloor$, where $0 < x < 1$, then the above expression is approximately
\begin{equation} \label{Brownfunc}
u(x) = E \bigg[ \exp \bigg( - \frac{A^{2(1 - 2^{-(j-1)})}}{2} \int_0^U \frac{1}{1 - B(s)} \: ds \bigg) {\bf 1}_{\{B(U) = 0\}} \bigg],
\end{equation}
where $U = \inf\{t: B(t) \in \{0, 1\} \}$ and $(B(t), t \geq 0)$ is Brownian motion started at $x$.  Here we are also using that $N^2 \gamma \sim A^{2(1 - 2^{-(j-1)})}$, where the factor of $N^2$ comes from the time change in replacing random walk with Brownian motion.  Since the probability that we get either fixation or a successful mutation is $1 - u(x)$, and we need to take a limit as the number of mutants at time zero gets small, we have $$\lambda_j = \lim_{x \rightarrow 0} \frac{1 - u(x)}{x}.$$  Thus, the problem reduces to evaluating the Brownian functional (\ref{Brownfunc}).  One can obtain a differential equation for $u(x)$ using the Feynman-Kac formula, and then get a series solution to the differential equation, from which the formula (\ref{lamj}) follows.  Details of this argument occupy section 6 of \cite{DSS}.

\subsection{Rapid mutations}

It remains to handle the case when $N \mu \nrightarrow 0$.  With this scaling, fixation will not occur before time $\tau_m$.  However, the waiting time between the type 1 mutation that will eventually produce a type $m$ descendant and the actual appearance of the type $m$ descendant can no longer be ignored.  As a result, waiting times are no longer sums of exponential random variables.  Instead, we obtain the following result.  The $m=2$ case of part 3 is equivalent to the special case of Theorem 1 in \cite{DSS} when $u_1 = u_2 = \mu$.

\begin{Theo}\label{bigmu}
We have the following limiting results when $N \mu \nrightarrow 0$.
\begin{enumerate}
\item If $\mu \gg N^{-2/m}$, then $$\lim_{N \rightarrow \infty} P(\tau_m > N^{-1/m} \mu^{-1} t) = \exp \bigg( - \frac{t^m}{m!} \bigg).$$

\item If $N^{-1/(1 + (m-j-2)2^{-(j+1)})} \ll \mu \ll N^{-1/(1 + (m-j-1)2^{-j})}$ for some $j = 1, \dots, m-2$, then $$\lim_{N \rightarrow \infty} P(\tau_m > N^{-1/(m-j)} \mu^{-1-(1 - 2^{-j})/(m-j)} t) = \exp \bigg( - \frac{t^{m-j}}{(m-j)!} \bigg).$$

\item If $\mu \sim A N^{-1/(1 + (m - j - 1)2^{-j})}$ for some $j = 1, \dots, m-1$ and some constant $A > 0$, then $$\lim_{N \rightarrow \infty} P(\tau_m > \mu^{-(1 - 2^{-j})} t) = \exp \bigg(- \frac{A^{1+(m-j-1)2^{-j}}}{(m-j-1)!} \int_0^t (t-s)^{m-j-1} \frac{1-e^{-2s}}{1+e^{-2s}} \: ds \bigg).$$
\end{enumerate}
\end{Theo}

We now explain the intuition behind these results.  Recall that $X_j(t)$ is the number of individuals with $j$ mutations at time $t$.  Because there are $N$ individuals getting mutations at rate $\mu$, we have $E[X_1(t)] \approx N \mu t$ for small $t$.  Each of these individuals acquires a second mutation at rate $\mu$, so $$E[X_2(t)] \approx \mu \int_0^t N \mu s \: ds = \frac{N \mu^2 t^2}{2}.$$ Repeating this reasoning, we get $E[X_j(t)] \approx N \mu^j t^j/j!$.

When the mutation rate is sufficiently large, there is a Law of Large Numbers, and the fluctuations in the number of individuals with $j$ mutations are small relative to $E[X_j(t)]$.  In this case, $X_j(t)$ is well approximated by its expectation.  When the mutation rate is sufficiently small, most of the time there are no individuals with $j$ mutations in the population, and when an individual gets a $j$th mutation, this mutation either dies out or, with probability $q_{m-j+1}$, produces a type $m$ descendant on a time scale much faster than $\tau_m$.  In this case, the problem reduces to determining how long we have to wait for a $j$th mutation that is successful enough to produce a type $m$ descendant.  There is also a borderline case in which we get stochastic effects in the limit both from the number of type $j$ individuals in the population and from the time between the appearance of a type $j$ individual that will eventually have a type $m$ descendant and the birth of the type $m$ descendant.

If the mutation rate is fast enough so that $X_{m-1}(t) \approx E[X_{m-1}(t)]$ up to time $\tau_m$, then since each individual with $m-1$ mutations gets an $m$th mutation at rate $\mu$, we get
\begin{equation} \label{biggesteq}
P(\tau_m > t) \approx \exp \bigg(- \mu \int_0^t \frac{N \mu^{m-1} s^{m-1}}{(m-1)!} \: ds \bigg) = \exp \bigg( -  \frac{N \mu^m t^m}{m!} \bigg).
\end{equation}
This leads to the result in part 1 of Theorem \ref{bigmu} if we substitute $N^{-1/m} \mu^{-1} t$ in place of $t$ in (\ref{biggesteq}).  In this regime, mutations are happening fast enough that births and deaths do not affect the limiting result, and we get the same result that we would get if $\tau_m$ were simply the first time that one of $N$ independent rate $\mu$ Poisson processes reaches the value $m$.  Consequently, as can be seen by integrating (\ref{adrate}), this result agrees with the result of Armitage and Doll \cite{armdoll}, who did not consider cell division and cell death in their original model.  The result when $m = 2$ agrees with a result in section 4 of \cite{imkn05}, and with (12.18) of \cite{nowak}.

Next, suppose mutation rates are fast enough so that $X_{m - j - 1}(t) \approx E[X_{m-j-1}(t)]$ up to time $\tau_m$, but slow enough that the time between the appearance of a ``successful" type $m - j$ individual that will have a type $m$ descendant and the birth of the type $m$ descendant is small relative to $\tau_m$.  Then each type $m - j - 1$ individual experiences ``successful" mutations at rate $\mu q_{j+1} \approx \mu^{2 - 2^{-j}}$ by Proposition \ref{family}, so
$$P(\tau_m > t) \approx \exp \bigg( -\mu^{2 - 2^{-j}} \int_0^t \frac{N \mu^{m-j-1} s^{m-j-1}}{(m-j-1)!} \: ds \bigg) = \exp \bigg( - \frac{N \mu^{m-j+1-2^{-j}} t^{m-j}}{(m-j)!} \bigg).$$  This leads to the result in part 2 of Theorem \ref{bigmu}.  The borderline cases are handled by part 3 of Theorem \ref{bigmu}.

To understand where the boundaries between the different types of behavior occur, first recall that
the number of type $k$ individuals born by time $t$ is of the order $N \mu^{k} t^{k}$.  Because each individual gives birth and dies at approximately rate one, the number of births and deaths of type $k$ individuals by time $t$ is of order $N \mu^{k} t^{k+1}$.  Because the standard deviation of the position of a random walk after $M$ steps is of order $M^{1/2}$, the standard deviation of the number of type $k$ individuals by time $t$ is of order $N^{1/2} \mu^{k/2} t^{(k+1)/2}$.  Therefore, we have $X_k(t) \approx E[X_k(t)]$ whenever $N^{1/2} \mu^{k/2} t^{(k+1)/2} \ll N \mu^k t^k$ or, equivalently, whenever $1 \ll N \mu^k t^{k-1}$.  See Proposition \ref{detprop} below for a precise statement of this result.

Each type $k$ individual experiences a mutation that will have a type $m$ descendant at rate $\mu q_{m-k} \approx \mu^{2 - 2^{-(m-k-1)}}$.  Therefore, the expected number of such mutations by time $t$ is of the order $N \mu^{k} t^{k} \cdot \mu^{2 - 2^{-(m-k-1)}} \cdot t = N \mu^{k+2-2^{-(m-k-1)}} t^{k+1}$.  This expression is of order one when $t$ is of order $N^{-1/(k+1)} \mu^{-1 - (1-2^{-(m-k-1)})/(k+1)}$, which is consequently the order of magnitude of the time we have to wait for one such mutation to occur.  It now follows from the result of the previous paragraph that $X_k(t) \approx E[X_k(t)]$ up to time $\tau_m$ whenever 
\begin{equation}\label{messO}
1 \ll N \mu^{k} (N^{-1/(k+1)} \mu^{-1 - (1-2^{-(m-k-1)})/(k+1)})^{k-1}.
\end{equation}
The expression on the right-hand side of (\ref{messO}) can be simplified to $(N^2 \mu^{2 + (k-1) 2^{-(m-k-1)}})^{1/(k+1)}$, so (\ref{messO}) is equivalent to the condition
\begin{equation}\label{heurmu}
\mu \gg N^{-1/(1 + (k-1)2^{-(m-k)})}.
\end{equation}
This condition can be compared to the condition for part 2 of Theorem \ref{bigmu}, which entails that (\ref{heurmu}) holds for $k = m-j-1$ but not for $k = m-j$, and therefore the number of type $m-j-1$ individuals, but not the number of type $m-j$ individuals, is approximately deterministic through time $\tau_m$.

If instead $\mu$ is of the order $N^{-1/(1 + (m-j-1)2^{-j})}$ for some $j = 1, \dots, m-1$, then on the relevant time scale the number of individuals of type $m-j-1$ behaves deterministically, but the number of individuals of type $m-j$ has fluctuations of the same order as the expected value.  As a result, there are stochastic effects from the number of type $m-j$ individuals in the population.  In this case, there are also stochastic effects from the time between the birth of type $m-j$ individual that will have a type $m$ descendant and the time that the type $m$ descendant is born.  Calculating the form of the limiting distribution in these borderline cases involves working with a two-type branching process.  This branching process is very similar to a process analyzed in chapter 3 of \cite{wk}, which explains the resemblance between part 3 of Theorem \ref{bigmu} and (3.20) of \cite{wk}.  Similar analysis using generating functions of branching processes that arise in multi-stage models of cancer has been carried out in \cite{mdv88, ml90, ml92}.  The work in \cite{ml90} allows for time-dependent parameters, while a three-stage model is analyzed in \cite{ml92}.

\subsection{The case $m = 3$}

To help the reader understand the different limiting behaviors, we summarize here the results when $m = 3$.  There are 9 different limiting regimes in this case; in general for the waiting time to get $m$ mutations, there are $4m-3$ limiting regimes.  Below $Z_1$ and $Z_2$ have the exponential distribution with mean one, and $Y_1$ and $Y_2$ have the exponential distributions with mean $\lambda_2$ and $\lambda_3$ respectively, where $\lambda_2$ and $\lambda_3$ are given by (\ref{lamj}).  The random variables $Z_1$, $Z_2$, $Y_1$, and $Y_2$ are assumed to be independent.

\begin{itemize}
\item If $\mu \ll N^{-2}$, then by part 1 of Theorem \ref{smallmu}, $\mu \tau_3 \rightarrow_d Z_1 + Z_2$.  We wait for two fixations, and then the third mutation happens quickly.

\item If $\mu \sim A N^{-2}$, then by the $j=2$ case of Theorem \ref{border}, $\mu \tau_3 \rightarrow_d Z_1 + Y_1$.  We wait for one fixation, then either a second fixation (after which the third mutation would happen quickly) or a second mutation that will not fixate but will have a descendant that gets a third mutation.

\item If $N^{-2} \ll \mu \ll N^{-4/3}$, then by the $j=2$ case of part 2 of Theorem \ref{smallmu}, $\mu \tau_3 \rightarrow_d Z_1$.  We wait for one fixation, and then the other two mutations happen quickly.

\item If $\mu \sim A N^{-4/3}$, then by the $j=3$ case of Theorem \ref{border}, $\mu \tau_3 \rightarrow_d Y_2$.  We wait either for a fixation (after which the other two mutations would happen quickly) or a mutation that will not fixate but will have a descendant with two additional mutations.

\item If $N^{-4/3} \ll \mu \ll N^{-1}$, then by part 3 of Theorem \ref{smallmu}, $N \mu^{7/4} \tau_3 \rightarrow_d Z_1$.  Fixation does not happen before time $\tau_3$, but we wait an exponentially distributed time for a mutation that is successful enough to have a descendant with three mutations.

\item If $\mu \sim A N^{-1}$, then by the $j=2$ case of part 3 of Theorem \ref{bigmu}, $$P(\mu^{3/4} \tau_3 > t) \rightarrow \exp \bigg( -A \int_0^t \frac{1 - e^{-2s}}{1 + e^{-2s}} \: ds \bigg).$$

\item If $N^{-1} \ll \mu \ll N^{-2/3}$, then by the $j=1$ case of part 2 of Theorem \ref{bigmu}, $P(N^{1/2} \mu^{5/4} \tau_3 > t) \rightarrow \exp(-t^2/2)$.  The number of individuals with one mutation is approximately deterministic, and the stochastic effect comes from waiting for a second mutation that is successful enough to have a descendant with a third mutation.

\item If $\mu \sim A N^{-2/3}$, then by the $j=1$ case of part 3 of Theorem \ref{bigmu}, $$P(\mu^{1/2} \tau_3 > t) \rightarrow \exp \bigg( - A^{3/2} \int_0^t (t-s) \frac{1 - e^{-2s}}{1 + e^{-2s}} \: ds \bigg).$$

\item If $\mu \gg N^{-2/3}$, then by part 1 of Theorem \ref{bigmu}, $P(N^{1/3} \mu \tau_3 > t) \rightarrow \exp(-t^3/6)$.  The number of individuals with two mutations is approximately deterministic, and the stochastic effect comes from waiting for the third mutation.
\end{itemize}

\subsection{Power law asymptotics and implications for cancer modeling} \label{powersec}

Because the probability that an individual develops a particular type of cancer during his or her lifetime is small, it seems unlikely that it will be possible to observe the full limiting distribution of the waiting time for cancer from data on cancer incidence.  Instead, we will observe only the left tail of this distribution.  Consequently, what is likely to be most relevant for applications are asymptotic formulas as $t \rightarrow 0$.  Throughout this subsection, write $f(t) \approx g(t)$ to mean that $f(t)/g(t) \rightarrow 1$ as $t \rightarrow 0$.  Recall that if $S_j$ is the sum of $j$ independent exponential random variables with mean one, then $P(S_j \leq t) \approx t^j/j!$.  This fact, combined with the approximation $1 - \exp(-t^{m-j}/(m-j)!) \approx t^{m-j}/(m-j)!$, allows us to deduce the following corollary of Theorems \ref{smallmu} and \ref{bigmu}.

\begin{Cor}\label{powercor}
We have the following asymptotic formulas as $t \rightarrow 0$:
\begin{enumerate}
\item If $\mu \ll N^{-2}$, then $$\lim_{N \rightarrow \infty} P(\tau_m \leq \mu^{-1} t) \approx \frac{t^{m-1}}{(m-1)!}.$$

\item If $N^{-2^{j-1}/(2^{j-1}-1)} \ll \mu \ll N^{-2^j/(2^j - 1)}$ for some $j = 2, \dots, m-1$, then 
$$P(\tau_m \leq \mu^{-1} t) \approx \frac{t^{m-j}}{(m-j)!}.$$

\item If $N^{-2^{m-1}/(2^{m-1}-1)} \ll \mu \ll N^{-1}$, then 
$$P(\tau_m \leq N^{-1} \mu^{-2 + 2^{-(m-1)}} t) \approx t.$$

\item If $N^{-1/(1 + (m-j-2)2^{-(j+1)})} \ll \mu \ll N^{-1/(1 + (m-j-1)2^{-j})}$ for some $j = 1, \dots, m-2$, then $$\lim_{N \rightarrow \infty} P(\tau_m \leq N^{-1/(m-j)} \mu^{-1-(1 - 2^{-j})/(m-j)} t) \approx \frac{t^{m-j}}{(m-j)!}.$$

\item If $\mu \gg N^{-2/m}$, then $$\lim_{N \rightarrow \infty} P(\tau_m \leq N^{-1/m} \mu^{-1} t) \approx
\frac{t^m}{m!}.$$
\end{enumerate}
\end{Cor}
By integrating (\ref{adrate}), we see that the result in part 5 of the corollary, which says that the probability of getting cancer by time $t$ behaves like $C t^m$, agrees with the result of Armitage and Doll.  However, parts 1 through 4 of the corollary show that in an $m$-stage model of cancer, the probability of getting cancer by time $t$ could behave like $C t^j$ for any $j = 1, 2, \dots, m$, depending on the relationship between $\mu$ and $N$.  This range of behavior can occur because not all of the $m$ events required for cancer are necessarily ``rate limiting".  For example, when part 2 of the corollary applies, there are $m-j$ fixations, and then the remaining $j$ mutations happen on a much faster time scale.  Consequently, it is not possible to deduce the number of mutations required for cancer just from the power law relationship between age and cancer incidence.

Corollary \ref{powercor} also shows that in our $m$-stage model, the probability of getting cancer by time $t$ will never behave like $C t^j$ for $j > m$.  However, as noted by Armitage and Doll (see \cite{arm85, armdoll}), higher powers could arise if the mutation rate, instead of being constant, increases over time like a power of $t$.  Also, the probability of getting cancer by time $t$ could increase more rapidly than $t^m$ if cells with mutations have a selective advantage over other cells, allowing their number to increase more rapidly than our model predicts.  This explains, in part, the success of two-stage models in fitting a wide variety of cancer incidence data, as documented in \cite{mknud81}.

\section{Proof of Theorems \ref{smallmu} and \ref{border}}\label{smallpf}

Recall that part 3 of Theorem \ref{smallmu} is a special case of Theorem 2 of \cite{DSS}, so we need to prove only parts 1 and 2.  We begin by recording three lemmas.  Lemma \ref{timelem}, which just restates (3.6), (3.8), and Lemma 3.1 of \cite{DSS}, bounds the amount of time that a mutation is in the population before it dies out or fixates.  Lemma \ref{qmprime} complements Proposition \ref{family}.  Lemma \ref{fastj} is a direct consequence of part 3 of Theorem \ref{smallmu}.  In these lemmas and throughout the rest of the paper, $C$ denotes a positive constant not depending on $N$ whose value may change from line to line.

\begin{Lemma}\label{timelem}
Consider a model of a population of size $N$ in which all individuals are either type 0 or type 1.  The population starts with just one type 1 individual and evolves according to the Moran model, so each individual dies at rate one and then gets replaced by a randomly chosen individual from the population.  Let $X(t)$ be the number of type 1 individuals at time $t$.  Let $T = \inf\{t: X(t) \in \{0, N\}\}$.  Let $L_k$ be the Lebesgue measure of $\{t: X(t) = k\}$.  Then for $k = 1, \dots, N-1$,
\begin{equation}\label{ELk}
E[L_k] = \frac{1}{k}.
\end{equation}
Also,
\begin{equation}\label{ET}
E[T] \leq C \log N
\end{equation}
and for all $0 \leq t \leq N$,
\begin{equation}\label{KolmPt}
P(T > t) \leq C/t.
\end{equation}
\end{Lemma}

\begin{Lemma}\label{qmprime}
Consider the model of Proposition \ref{family}.  Let $q_m'$ be the probability that a type $m$ individual is born at some time, but that eventually all individuals have type zero.  Suppose $N \mu^{1-2^{-(m-1)}} \rightarrow 0$ as $N \rightarrow \infty$.  Then $$q_m' \ll 1/N.$$
\end{Lemma}

\begin{proof}
The event that all individuals eventually have type zero has probability $(N-1)/N$ regardless of the mutation rate.  On this event, reducing the mutation rate can only reduce the probability of eventually getting a type $m$ individual.  Therefore, it suffices to prove the result when
\begin{equation}\label{muassum}
N \mu^{1 - 2^{-(m-2)}} \rightarrow \infty.
\end{equation}
If a type $m$ individual eventually is born, then some type 2 mutation must have a type $m$ descendant.  By (\ref{ELk}), for $k = 1, \dots, N-1$, the expected amount of time for which there are $k$ individuals of nonzero type is $1/k$.  While there are $k$ individuals of nonzero type, type 2 mutations occur at rate at most $k \mu$.  On the event that there is no fixation, the number of individuals of nonzero type never reaches $N$, and the expected number of type 2 mutations while there are fewer than $N$ individuals of nonzero type is at most
$$\sum_{k=1}^{N-1} \frac{1}{k} \cdot k \mu \leq N \mu.$$  When (\ref{muassum}) holds, we can apply Proposition \ref{family} to see that if $m \geq 3$ then each type 2 mutation has probability at most $C \mu^{1 - 2^{-(m-2)}}$ of having a type $m$ descendant.  This inequality holds trivially if $m = 2$.  It follows that $$q_m' \leq (N \mu) (C \mu^{1 - 2^{-(m-2)}}) = C N \mu^{2 - 2^{-(m-2)}},$$
and therefore $N q_m' \leq C(N \mu^{1 - 2^{-(m-1)}})^2 \rightarrow 0$, as claimed.
\end{proof}

\begin{Lemma}\label{fastj}
Suppose $j \geq 2$.  If $N^{-2^{j-1}/(2^{j-1} - 1)} \ll \mu \ll 1/N$, then for all $\epsilon > 0$,
$$\lim_{N \rightarrow \infty} P(\tau_j < \epsilon \mu^{-1}) = 1.$$
\end{Lemma}

\begin{proof}
Part 3 of Theorem \ref{smallmu} gives $\lim_{N \rightarrow \infty} P(N \mu^{2 - 2^{-(j-1)}} \tau_j \leq t) = 1 - e^{-t}$ for all $t > 0$.  The result follows immediately because $\mu \ll N \mu^{2 - 2^{-(j-1)}}$ by assumption.
\end{proof}

\begin{proof}[Proof of parts 1 and 2 of Theorem \ref{smallmu}]
Suppose either $j = 1$ and $\mu \ll N^{-2}$, or $j = 2, \dots, m-1$ and $N^{-2^{j-1}/(2^{j-1}-1)} \ll \mu \ll N^{-2^j/(2^j-1)}$.   Let $\gamma_i$ be the time of the $i$th mutation, so the points $(\gamma_i)_{i=1}^{\infty}$ form a rate $N \mu$ Poisson process on $[0, \infty)$.  Call the $i$th mutation bad if at time $\gamma_i$, there is another mutation in the population that has not yet died out or fixated.  Otherwise, call the mutation good.  For all $i$, let $\xi_i = 1$ if the $i$th mutation fixates, and let $\xi_i = 0$ otherwise.  We have $P(\xi_i = 1) = 1/N$ for all $i$, but the random variables $(\xi_i)_{i=1}^{\infty}$ are not independent because if two mutations are present at the same time on different individuals, at most one of the mutations can fixate.

Let $({\tilde \xi}_i)_{i=1}^{\infty}$ be a sequence of i.i.d. random variables, independent of the population process, such that $P({\tilde \xi}_i = 1) = 1/N$ and $P({\tilde \xi}_i = 0) = (N-1)/N$ for all $i$.  Define another sequence $(\xi_i')_{i=1}^{\infty}$ such that $\xi_i' = \xi_i$ if the $i$th mutation is good and $\xi_i' = {\tilde \xi}_i$ if the $i$th mutation is bad.  If the $i$th mutation is good, then $P(\xi_i = 1|(\xi'_k)_{k=1}^{i-1}) = 1/N$, so $(\xi_i')_{i=1}^{\infty}$ is an i.i.d. sequence.  Let $\sigma_1 = \inf\{\gamma_i: \xi_i = 1\}$ and for $k \geq 2$, let $\sigma_k = \inf\{\gamma_i > \sigma_{k-1}: \xi_i = 1\}$.  Likewise, let $\sigma_1' = \inf\{\gamma_i: \xi_i' = 1\}$ and for $k \geq 2$, let $\sigma_k' = \inf\{\gamma_i > \sigma_{k-1}: \xi_i' = 1\}$.  The points $\gamma_i$ for which $\xi_i' = 1$ form a Poisson process of rate $\mu$, so $\mu \sigma_{m-j}'$ has the gamma distribution with parameters $(m-j, 1)$.

Let $\epsilon > 0$, and choose $t$ large enough that 
\begin{equation}\label{sigmuep}
P(\sigma_{m-j}' > \mu^{-1} t) < \epsilon.
\end{equation}
Note that because $\mu \sigma_{m-j}'$ has a gamma distribution for all $N$, here $t$ does not depend on $N$.  The expected number of mutations by time $\mu^{-1} t$ is $(N \mu)(\mu^{-1} t) = Nt$.  After a mutation occurs, the number of individuals descended from this mutant individual evolves in the same way as the number of type 1 individuals in Lemma \ref{timelem}.  Therefore, by (\ref{ET}), the expected amount of time, before time $\mu^{-1}t$, that there is a mutation in the population that has not yet disappeared or fixated is at most $C (N \log N) t$.  Therefore, the expected number of bad mutations before time $\mu^{-1} t$ is at most $(N \mu)(C (N \log N) t) = C (N^2 \log N) \mu t$.  If a bad mutation occurs at time $\gamma_i$, the probability that either $\xi_i$ or $\xi_i'$ equals one is at most $2/N$, so
$$P(\xi_i = \xi_i' \mbox{ for all }i \mbox{ such that }\gamma_i \leq \mu^{-1} t) \geq 1 - 2C (N \log N) \mu t.$$  Because $\mu \ll 1/(N \log N)$, it follows by letting $\epsilon \rightarrow 0$ that
\begin{equation}\label{sigpr}
\lim_{N \rightarrow \infty} P(\sigma_{m-j}' = \sigma_{m-j}) =1.
\end{equation}
Thus, $\mu \sigma_{m-j} \rightarrow_d S_{m-j}.$  To complete the proof, it remains to show that
\begin{equation}\label{taumsmall}
\mu (\tau_m - \sigma_{m-j}) \rightarrow_p 0.
\end{equation}

We first prove that
\begin{equation}\label{taumbig}
\lim_{N \rightarrow \infty} P(\tau_m < \sigma_{m-j}) = 0
\end{equation}
If $\tau_m < \sigma_{m-j}$, then before time $\sigma_{m-j}$, there must be a type $k$ mutation for some $k \leq m-j$ that does not fixate but has a type $m$ descendant.  We will bound the probability of this event.  Recall that the expected number of mutations before time $\mu^{-1}t$ is $Nt$.  Because $\mu \ll N^{-2^j/(2^j - 1)}$, we can apply Lemma \ref{qmprime} with $j+1$ in place of $m$ to get that the probability that a type $m-j$ mutation does not fixate but has a type $m$ descendant is asymptotically much smaller than $1/N$.  Thus, the probability that before time $\mu^{-1}t$, there is a type $k$ mutation for some $k \leq m-j$ that does not fixate but has a type $m$ descendant is asymptotically much smaller than $(N t)(1/N)$, and therefore goes to zero as $N \rightarrow \infty$.  Combining this result with (\ref{sigmuep}) and (\ref{sigpr}) gives (\ref{taumbig}).

We now prove (\ref{taumsmall}).  Choose $\epsilon > 0$.  Let ${\tilde \gamma}_i$ be the time when the mutation at time $\gamma_i$ disappears or fixates.  By (\ref{ET}), we have $E[{\tilde \gamma}_i - \gamma_i] \leq C \log N$.  It follows from Markov's Inequality that $P({\tilde \gamma}_i - \gamma_i > \mu^{-1} \epsilon) \leq C \log N/(\mu^{-1} \epsilon)$.  Because the expected number of mutations by time $\mu^{-1}t$ is $Nt$, another application of Markov's Inequality gives $$P({\tilde \gamma}_i - \gamma_i > \mu^{-1} \epsilon \mbox{ for some }i \mbox{ such that }\gamma_i < \mu^{-1} t) \leq N t \cdot \frac{C \log N}{\mu^{-1} \epsilon} = \frac{Ct}{\epsilon}(N \log N) \mu,$$ which goes to zero as $N \rightarrow \infty$.  Therefore, in view of (\ref{sigmuep}) and (\ref{sigpr}), if $\zeta$ is the time when the mutation at time $\sigma_{m-j}$ fixates, we have
\begin{equation}\label{muzetsig}
\mu(\zeta - \sigma_{m-j}) \rightarrow_p 0
\end{equation}
Now (\ref{taumsmall}) will be immediate from (\ref{taumbig}) and (\ref{muzetsig}) once we show that for all $\epsilon > 0$,
\begin{equation}\label{zetaep}
\lim_{N \rightarrow \infty} P(\mu (\tau_m - \zeta) > \epsilon) = 0.
\end{equation}
When $j \geq 2$, equation (\ref{zetaep}) follows from Lemma \ref{fastj} because after time $\sigma_{m-j}$, at most $j$ more mutations are needed before we reach time $\tau_m$.  When $j = 1$, we reach the time $\tau_m$ as soon as there is another mutation after time $\sigma_{m-j}$, so $\tau_m - \zeta$ is stochastically dominated by an exponentially distributed random variable with rate $N \mu$.  It follows that (\ref{zetaep}) holds in this case as well.
\end{proof}

Most of the work involved in proving Theorem \ref{border} is contained in the proof of the following result, which is a special case of Lemma 7.1 of \cite{DSS}.

\begin{Lemma}
Suppose $\mu \sim A N^{-2^{j-1}/(2^{j-1} - 1)}$ for some $j = 2, \dots, m$ and some constant $A > 0$.
Consider the model of Proposition \ref{family}.  Let $q_j'$ be the probability that either a type $j$ individual is born at some time, or eventually all individuals in the population have type greater than zero.  Then $\lim_{N \rightarrow \infty} N q_j' = \lambda_j$, where $\lambda_j > 1$ is given by (\ref{lamj}).
\end{Lemma}

\begin{proof}[Proof of Theorem \ref{border}]
The proof is similar to the proof of parts 1 and 2 of Theorem \ref{smallmu}.  Define the sequences $(\gamma_i)_{i=1}^{\infty}$, $(\xi_i)_{i=1}^{\infty}$, $({\tilde \xi}_i)_{i=1}^{\infty}$ and $(\xi_i')_{i=1}^{\infty}$ as in the proof of parts 1 and 2 of Theorem \ref{smallmu}.  Also define a sequence $(\zeta_i)_{i=1}^{\infty}$ of $\{0, 1\}$-valued random variables such that $\zeta_1 = 1$ if the mutation at time $\gamma_i$ either fixates or has a descendant that gets $j-1$ additional mutations.  Let $({\tilde \zeta}_i)_{i=1}^{\infty}$ be a sequence of i.i.d. random variables, independent of the population process, such that $P({\tilde \zeta}_i = 1) = \lambda_j/N$ and $P({\tilde \zeta}_i = 0) = (N - \lambda_j)/N$ for all $i$, and ${\tilde \zeta}_i = 1$ whenever ${\tilde \xi}_i = 1$.  Let $\zeta_i' = \zeta_i$ if the $i$th mutation is good, and let $\zeta_i' = {\tilde \zeta}_i$ otherwise.  Let $\sigma_0 = 0$.  For $k = 1, \dots, m-j$, let $\sigma_k = \inf\{\gamma_i > \sigma_{k-1}: \xi_i = 1\}$.  Let $\sigma_{m-j+1} = \inf\{\gamma_i > \sigma_{m-j}: \zeta_i = 1\}$.  Define $\sigma_1', \dots, \sigma_{m-j+1}'$ in the same way using the random variables $\xi_i'$ and $\zeta_i'$.  It is clear from the construction that $\sigma_{m-j+1}'$ has the same distribution as $S_{m-j} + Y$.  By the same argument used in the proof of parts 1 and 2 of Theorem \ref{smallmu}, with a bound of $2 \lambda_j/N$ replacing the bound of $2/N$, we get $$\lim_{N \rightarrow \infty} P(\sigma'_{m-j+1} = \sigma_{m-j+1}) = 1,$$ which implies $\mu \sigma_{m-j+1} \rightarrow_d S_{m-j} + Y$.  This argument also gives that the mutation at time $\sigma_{m-j+1}$ is good with probability tending to one as $N \rightarrow \infty$.

We next claim that
\begin{equation}\label{taumbig2}
\lim_{N \rightarrow \infty} P(\tau_m < \sigma_{m-j+1}) = 0.
\end{equation}
If $\sigma_{m-j} < \gamma_i < \sigma_{m-j+1}$, then by the definition of $\sigma_{m-j+1}$, no descendant of the mutation at time $\gamma_i$ can have a type $m$ descendant.
Therefore, if $\tau_m < \sigma_{m-j+1}$, then before time $\sigma_{m-j}$ there must be a type $k$ mutation for some $k \leq m-j$ that does not fixate but has a type $m$ descendant.  Because $\mu \ll N^{-2^j/(2^j-1)}$, the probability of this event goes to zero by the same argument given in the proof of parts 1 and 2 of Theorem \ref{smallmu}, which implies (\ref{taumbig2}).

It remains only to prove
\begin{equation}\label{taumsmall2}
\mu(\tau_m - \sigma_{m-j+1}) \rightarrow_p 0.
\end{equation}
Let $\epsilon > 0$, and choose $t$ large enough that $P(\sigma'_{m-j+1} > \mu^{-1} t) < \epsilon$.  Let $\epsilon > 0$.  By the same argument given in the proof of parts 1 and 2 of Proposition \ref{smallmu}, the probability that some mutation before time $\mu^{-1}t$ takes longer than $\mu^{-1} \epsilon$ to die out or fixate tends to zero as $N \rightarrow \infty$.  Therefore, if $\zeta$ is the time when the mutation at time $\sigma_{m-j+1}$ dies out or fixates, then $\mu(\zeta - \sigma_{m-j+1}) \rightarrow_p 0$.  If the mutation at time $\sigma_{m-j+1}$ fixates, then only $j-1$ more mutations are needed before we reach time $\tau_m$.  Therefore, conditional on this fixation, when $j \geq 3$ we get $\mu(\tau_m - \zeta) \rightarrow_p 0$ by applying Lemma \ref{fastj} with $j-1$ in place of $j$, while the result $\mu(\tau_m - \zeta) \rightarrow_p 0$ is immediate when $j = 2$.  Alternatively, if the mutation at time $\sigma_{m-j+1}$ does not fixate and the mutation at time $\sigma_{m-j+1}$ is good, then $\tau_m \leq \zeta$.  Because the mutation at time $\sigma_{m-j+1}$ is good with probability tending to one as $n \rightarrow \infty$, we conclude (\ref{taumsmall2}).
\end{proof}

\section{Proof of parts 1 and 2 of Theorem \ref{bigmu}}

The first step in the proof of Theorem \ref{bigmu} is to establish conditions, stated in Proposition \ref{detprop} below, under which the number of type $k$ individuals is essentially deterministic, in the sense that it can be well approximated by its expectation.  It will follow that when $\mu \gg N^{-2/m}$, the number of individuals with type $m-1$ is approximately deterministic until time $\tau_m$.  Since each type $m-1$ individual experiences a type $m$ mutation at rate $\mu$, the approximately deterministic behavior of the type $m-1$ individuals leads easily to a proof of part 1 of Theorem \ref{bigmu}.  When instead $N^{-1/(1 + (m-j-2)2^{-(j+1)})} \ll \mu \ll N^{-1/(1 + (m-j-1)2^{-j})}$, the number of individuals of type $m-j-1$ is approximately deterministic up to time $\tau_m$, as will be shown in Lemma \ref{GNlem} below.  The remainder of the proof of part 2 of Theorem \ref{bigmu} involves using a Poisson approximation technique to calculate the distribution of the time we have to wait for one of the type $m-j-1$ individuals to have a type $m-j$ mutation that will give rise to a type $m$ descendant.

We begin with a lemma bounding the expected number of type $k$ individuals.  Recall that $X_j(t)$ denotes the number of type $j$ individuals at time $t$, and $X_j(0) = 0$ for all $j \geq 1$.

\begin{Lemma}
Let $Y_k(t) = \sum_{j=k}^{\infty} X_j(t)$ be the number of individuals of type $k$ or higher at time $t$.
For all $k \geq 0$ and $t \geq 0$, we have $E[X_k(t)] \leq E[Y_k(t)] \leq N \mu^k t^k/k!$.
\label{expk}
\end{Lemma}

\begin{proof}
The first inequality is obvious, so it suffices to show $E[Y_k(t)] \leq N \mu^k t^k/k!$.  We proceed by induction.  Since $Y_0(t) \leq N$ for all $t \geq 0$, the result is true for $k = 0$.  Suppose $k \geq 1$ and $E[Y_{k-1}(t)] \leq N \mu^{k-1} t^{k-1}/(k-1)!$ for all $t \geq 0$.  The expected number of type $k$ mutations before time $t$ is at most $$\mu \int_0^t E[X_{k-1}(s)] \: ds \leq \int_0^t \frac{N \mu^k s^{k-1}}{(k-1)!} \: ds = \frac{N \mu^k t^k}{k!}.$$  Because individuals of type $k$ and higher give birth and die at the same rate, it follows that $E[Y_k(t)] \leq N \mu^k t^k/k!$.
\end{proof}

\begin{Prop} \label{detprop}
Suppose $k \geq 0$ and $T$ is a time that depends on $N$.  Assume that as $N \rightarrow \infty$, we have $\mu T \rightarrow 0$, $N \mu^k T^{k-1} \rightarrow \infty$, and $N \mu^k T^k \rightarrow \infty$.  Then for all $\epsilon > 0$,
\begin{equation}
\lim_{N \rightarrow \infty} P \bigg( \max_{0 \leq t \leq T} \bigg| X_k(t) - \frac{N \mu^k t^k}{k!} \bigg| > \epsilon N \mu^k T^k \bigg) = 0.
\label{deteq}
\end{equation}
\end{Prop}

\begin{proof}
We prove the result by induction and begin with $k = 0$.  Individuals of type one or higher are always being born and dying at the same rate.  Since new individuals of type one or higher also appear because of type 1 mutations, the process $(N - X_0(t), t \geq 0)$ is a bounded submartingale.  Let $\zeta = \inf\{t: N - X_0(t) > \epsilon N\}$.  By the Optional Sampling Theorem, we have $E[N - X_0(T)|\zeta \leq T] \geq \epsilon N$.  Since the rate of type 1 mutations is always bounded by $N \mu$, we have $E[N - X_0(T)] \leq N \mu T$.  Therefore, $$P \bigg( \max_{0 \leq t \leq T} |X_0(t) - N| > \epsilon N \bigg) = P(\zeta \leq T) \leq \frac{E[N - X_0(T)]}{E[N - X_0(T)|\zeta \leq T]} \leq \frac{N \mu T}{\epsilon N} \rightarrow 0$$ as $N \rightarrow \infty$ because $\mu T \rightarrow 0$.  It follows that when $k = 0$, (\ref{deteq}) holds for all $\epsilon > 0$.

Let $k \geq 1$.  Assume that (\ref{deteq}) holds with $k-1$ in place of $k$.  Let $B_k(t)$ be the number of type $k$ mutations up to time $t$.  Let $S_k(t)$ be the number of times, until time $t$, that a type $k$ individual gives birth minus the number of times that a type $k$ individual dies.  Note that $X_k(t) = B_k(t) - B_{k+1}(t) + S_k(t)$, so
\begin{equation}
\bigg| X_k(t) - \frac{N \mu^k t^k}{k!} \bigg| \leq B_{k+1}(t) + |S_k(t)| + \bigg| B_k(t) - \frac{N \mu^k t^k}{k!} \bigg|.
\label{3abstms}
\end{equation}
Therefore, it suffices to show that with probability tending to one as $N \rightarrow \infty$, the three terms on the right-hand side of (\ref{3abstms}) stay below $\epsilon N \mu^k T^k/3$ for $t \leq T$.

By Lemma \ref{expk}, for $0 \leq t \leq T$,
$$E[B_{k+1}(t)] = \mu \int_0^T E[X_k(t)] \: dt \leq \frac{N \mu^{k+1} T^{k+1}}{(k+1)!}.$$  By Markov's Inequality,
\begin{equation}
P\bigg( \max_{0 \leq t \leq T} B_{k+1}(t) > \frac{\epsilon}{3} N \mu^k T^k \bigg) = P \bigg( B_{k+1}(T) >  \frac{\epsilon}{3} N \mu^k T^k \bigg) \leq \frac{3 \mu T}{\epsilon (k+1)!} \rightarrow 0
\label{tm1to0}
\end{equation}
as $N \rightarrow \infty$ because $\mu T \rightarrow 0$.

Note that $S(0) = 0$, and since type $k$ individuals give birth and die at the same rate, the process $(S(t), 0 \leq t \leq T)$ is a martingale.  By Wald's Second Equation, $E[S(T)^2]$ is the expected number of births plus deaths of type $k$ individuals (not counting replacements of a type $k$ individual by another type $k$ individual) up to time $T$, which by Lemma \ref{expk} is at most
$$2 \int_0^T E[X_k(t)] \: dt \leq \frac{2 N \mu^k T^{k+1}}{(k+1)!}.$$  Therefore, by the $L^2$-Maximal Inequality for martingales, $$E \big[ \max_{0 \leq t \leq T} |S(t)|^2 \big] \leq 4 E[S(T)^2] \leq \frac{8 N \mu^k T^{k+1}}{(k+1)!}.$$  Now using Chebyshev's Inequality,
\begin{equation}
P \bigg( \max_{0 \leq t \leq T} |S_k(t)| > \frac{\epsilon}{3} N \mu^k T^k \bigg) \leq \frac{8 N \mu^k T^{k+1}}{(k+1)!} \bigg( \frac{3}{\epsilon N \mu^k T^k} \bigg)^2 = \frac{72}{(k+1)! N \mu^k T^{k-1}} \rightarrow 0
\label{tm2to0}
\end{equation}
as $N \rightarrow \infty$ because $N \mu^k T^{k-1} \rightarrow \infty$.

To bound the third term in (\ref{3abstms}), note that type $k-1$ individuals mutate to type $k$ at rate $\mu$.  Therefore, there exist inhomogeneous Poisson processes $(N_1(t), t \geq 0)$ and $(N_2(t), t \geq 0)$ whose intensities at time $t$ are given by $N \mu^k t^{k-1}/(k-1)! - \epsilon N \mu^k T^{k-1}/6$ and $N \mu^k t^{k-1}/(k-1)! + \epsilon N \mu^k T^{k-1}/6$ respectively such that on the event that
\begin{equation}
\max_{0 \leq t \leq T} \bigg| X_{k-1}(t) - \frac{N \mu^{k-1}t^{k-1}}{(k-1)!} \bigg| \leq \frac{\epsilon}{6} N \mu^{k-1} T^{k-1},
\label{indevent}
\end{equation}
we have $N_1(t) \leq B_k(t) \leq N_2(t)$ for $0 \leq t \leq T$.  To achieve this coupling, one can begin with points at the times of type $k$ mutations.  To get $(N_1(t), t \geq 0)$, when there is a type $k$ mutation at time $t$, remove this point with probability $[N \mu^k t^{k-1}/(k-1)! - \epsilon N \mu^k T^{k-1}/6]/\mu X_{k-1}(t-)$.  To get $(N_2(t), t \geq 0)$, add points of a time-inhomogeneous Poisson process whose rate at time $t$ is $[N \mu^k t^{k-1}/(k-1)! + \epsilon N \mu^k T^{k-1}/6] - \mu X_{k-1}(t)$.

Note that
\begin{equation}
E[N_1(t)] = \int_0^t \bigg( \frac{N \mu^k s^{k-1}}{(k-1)!} - \frac{\epsilon N \mu^k T^{k-1}}{6} \bigg) \: ds = \frac{N \mu^k t^k}{k!} - \frac{\epsilon}{6} N \mu^k T^{k-1} t
\label{EN1}
\end{equation}
and likewise
$$E[N_2(t)] = \frac{N \mu^k t^k}{k!} + \frac{\epsilon}{6} N \mu^k T^{k-1} t.$$
The process $(N_1(t) - E[N_1(t)], t \geq 0)$ is a martingale, and
\begin{equation}\label{MGvar}
E \big[ (N_1(T) - E[N_1(T)])^2 \big] = E[N_1(T)] = \frac{N \mu^k T^k}{k!} - \frac{\epsilon}{6} N \mu^k T^k.
\end{equation}
Therefore, Chebyshev's Inequality and the $L^2$-Maximal Inequality for martingales give
\begin{align}
P \bigg( \max_{0 \leq t \leq T} \big|N_1(t) - E[N_1(t)]\big| > \frac{\epsilon}{6} N \mu^k T^k \bigg) &\leq \frac{36 E \big[ \max_{0 \leq t \leq T} |N_1(t) - E[N_1(t)]|^2 \big]}{(\epsilon N \mu^k T^k)^2} \nonumber \\
&\leq \frac{144 E \big[(N_1(T) - E[N_1(T)])^2 \big]}{(\epsilon N \mu^k T^k)^2} \rightarrow 0
\label{VarN1}
\end{align}
as $N \rightarrow \infty$ by (\ref{MGvar}) because $N \mu^k T^k \rightarrow \infty$.  Combining (\ref{EN1}) with (\ref{VarN1}) gives
\begin{equation}\label{newN1}
\lim_{N \rightarrow \infty} P \bigg( \max_{0 \leq t \leq T} \bigg| N_1(t) - \frac{N \mu^k t^k}{k!} \bigg| > \frac{\epsilon}{3} N \mu^k T^k \bigg) = 0. 
\end{equation}
The same argument gives
\begin{equation}\label{newN2}
\lim_{N \rightarrow \infty} P \bigg( \max_{0 \leq t \leq T} \bigg| N_2(t) - \frac{N \mu^k t^k}{k!} \bigg| > \frac{\epsilon}{3} N \mu^k T^k \bigg) = 0. 
\end{equation}
as $N \rightarrow \infty$.  By the induction hypothesis, the event in (\ref{indevent}) occurs with probability tending to one as $N \rightarrow \infty$, so $N_1(t) \leq B_k(t) \leq N_2(t)$ for $0 \leq t \leq T$ with probability tending to one as $N \rightarrow \infty$.  Therefore, equations (\ref{newN1}) and (\ref{newN2}) imply that
\begin{equation}
\lim_{N \rightarrow \infty} P \bigg( \max_{0 \leq t \leq T} \bigg| B_k(t) - \frac{N \mu^k t^k}{k!} \bigg| > \frac{\epsilon}{3} N \mu^k T^k \bigg) = 0.
\label{tm3to0}
\end{equation}
The result follows from (\ref{3abstms}), (\ref{tm1to0}), (\ref{tm2to0}), and (\ref{tm3to0}).
\end{proof}

\begin{proof}[Proof of part 1 of Theorem \ref{bigmu}]
Suppose $\mu \gg N^{-2/m}$, and let $T = N^{-1/m} \mu^{-1} t$.  As $N \rightarrow \infty$, we have $\mu T = N^{-1/m} t \rightarrow 0$, $N \mu^{m-1} T^{m-2} = N^{2/m} \mu t^{m-2} \rightarrow \infty$, and $N \mu^{m-1} T^{m-1} = N^{1/m} t^{m-1} \rightarrow \infty$.  Therefore, by Proposition \ref{detprop}, if $\epsilon > 0$, then with probability tending to one as $N \rightarrow \infty$,
\begin{equation}\label{concevent}
\max_{0 \leq s \leq T} \bigg| X_{m-1}(s) - \frac{N \mu^{m-1} s^{m-1}}{(m-1)!} \bigg| \leq \epsilon N \mu^{m-1} T^{m-1}.
\end{equation}
Because each type $m-1$ individual experiences a type $m$ mutation at rate $\mu$, the random variable $$V = \int_0^{\tau_m} \mu X_{m-1}(s) \: ds$$ has an exponential distribution with mean one.  When (\ref{concevent}) holds, we have
$$\frac{N \mu^m T^m}{m!} - \epsilon N \mu^m T^m \leq \int_0^T \mu X_{m-1}(s) \: ds \leq \frac{N \mu^m T^m}{m!} + \epsilon N \mu^m T^m.$$  It follows that
\begin{align}
\limsup_{N \rightarrow \infty} P(\tau_m > T) &\leq \limsup_{N \rightarrow \infty} P \bigg( V > \frac{N \mu^m T^m}{m!} - \epsilon N \mu^m T^m \bigg) \nonumber \\
&= P \bigg( W > \frac{t^m}{m!} - \epsilon t^m \bigg) = \exp \bigg( - \frac{t^m}{m!} + \epsilon t^m \bigg), \nonumber
\end{align}
and likewise
$$\liminf_{N \rightarrow \infty} P(\tau_m > T) \geq \liminf_{N \rightarrow \infty} P \bigg( V > \frac{N \mu^m T^m}{m!} + \epsilon N \mu^m T^m \bigg) = \exp \bigg( - \frac{t^m}{m!} - \epsilon t^m \bigg).$$  Because these bounds hold for all $\epsilon > 0$, the result follows.
\end{proof}

We now work towards proving part 2 of Theorem \ref{bigmu}.  For the rest of this section, we assume that 
\begin{equation} \label{mucond}
N^{-1/(1 + (m-j-2)2^{-(j+1)})} \ll \mu \ll N^{-1/(1 + (m-j-1)2^{-j})}
\end{equation}
for some $j = 1, \dots, m-2$.  This condition implies that $N \mu \rightarrow \infty$ and $\mu \rightarrow 0$ as $N \rightarrow \infty$, and therefore
\begin{equation}\label{mujcond}
N \mu^{1-2^{-j}} \rightarrow \infty.
\end{equation}
Also, for the rest of this section, $t$ is fixed and
\begin{equation}\label{defT}
T = N^{-1/(m-j)} \mu^{-1-(1-2^{-j})/(m-j)} t.
\end{equation}
This means that
\begin{equation}\label{latenote}
N \mu^{m-j} T^{m-j} = \mu^{-(1-2^{-j})} t^{m-j}.
\end{equation}
Let $\epsilon > 0$.  Let $G_N$ be the event that $$\max_{0 \leq s \leq T} \bigg|X_{m-j-1}(s) - \frac{N \mu^{m-j-1} s^{m-j-1}}{(m-j-1)!} \bigg| \leq \epsilon N \mu^{m-j-1} T^{m-j-1}.$$  The next lemma shows that $G_N$ occurs with high probability, indicating that on the time scale of interest, the number of individuals with $m-j-1$ mutations stays close to its expectation. 

\begin{Lemma}\label{GNlem}
We have $\lim_{N \rightarrow \infty} P(G_N) = 1$.
\end{Lemma}

\begin{proof}
We need to verify the conditions of Proposition \ref{detprop} with $m-j-1$ in place of $k$.  By (\ref{mujcond}), as $N \rightarrow \infty$,
\begin{equation}\label{check1}
\mu T = N^{-1/(m-j)} \mu^{-(1-2^{-j})/(m-j)}t = (N \mu^{1-2^{-j}})^{-1/(m-j)}t \rightarrow 0.
\end{equation}
Also, using the first inequality in (\ref{mucond}),
\begin{align} \label{check2}
N \mu^{m-j-1} T^{m-j-2} &= N^{1 - (m-j-2)/(m-j)} \mu^{m-j-1-(m-j-2)-(m-j-2)(1-2^{-j})/(m-j)} t^{m-j-2} \nonumber \\
&= N^{2/(m-j)} \mu^{2/(m-j) + (m-j-2)2^{-j}/(m-j)} t^{m-j-2} \nonumber \\
&= (N \mu^{1 + (m-j-2)2^{-(j+1)}})^{2/(m-j)} t^{m-j-2} \rightarrow \infty.
\end{align}
Using the second inequality in (\ref{mucond}) and the fact that $m-j+1-2^{-j} > 1 + (m-j-1)2^{-j}$,
$$T = (N \mu^{m-j+1-2^{-j}})^{-1/(m-j)} t \gg (N^{1 - (m-j+1-2^{-j})/(1+(m-j-1)2^{-j})})^{-1/(m-j)} t \rightarrow \infty.$$
This result and (\ref{check2}) imply $N \mu^{m-j-1} T^{m-j-1} \rightarrow \infty$, which, in combination with (\ref{check1}) and (\ref{check2}), gives the lemma.
\end{proof}

The rest of the proof of part 2 of Theorem \ref{bigmu} is similar to the proof of Theorem 2 in \cite{DSS}.  It depends on the following result on Poisson approximation, which is part of Theorem 1 of \cite{agg} and was used also in \cite{DSS}.

\begin{Lemma}
Suppose $(A_i)_{i \in {\cal I}}$ is a collection of events, where ${\cal I}$ is any index set.  Let $W = \sum_{i \in {\cal I}} {\bf 1}_{A_i}$ be the number of events that occur, and let $\lambda = E[W] = \sum_{i \in {\cal I}} P(A_i)$.  Suppose for each $i \in {\cal I}$, we have $i \in \beta_i \subset {\cal I}$.  Let ${\cal F}_i = \sigma((A_j)_{j \in {\cal I} \setminus \beta_i})$.  Define
\begin{align}
b_1 &= \sum_{i \in {\cal I}} \sum_{j \in \beta_i} P(A_i) P(A_j), \nonumber \\
b_2 &= \sum_{i \in {\cal I}} \sum_{i \neq j \in \beta_i} P(A_i \cap A_j), \nonumber \\
b_3 &= \sum_{i \in {\cal I}} E \big[| P(A_i|{\cal F}_i) - P(A_i)|\big]. \nonumber
\end{align}
Then $|P(W = 0) - e^{-\lambda}| \leq b_1 + b_2 + b_3$.
\label{Poisson}
\end{Lemma}

We will use the next lemma to get the second moment estimate needed to bound $b_2$. When we apply this result, the individuals born at times $t_1$ and $t_2$ will both have the same type. We use different types in the statement of the lemma to make it easier to distinguish the descendants of the two individuals.  This result is Lemma 5.2 of \cite{DSS}.

\begin{Lemma}\label{2families}
Fix times $t_1 < t_2$.  Consider a population of size $N$ which evolves according to the Moran model in which all individuals initially have type 0.  There are no mutations, except that one individual becomes type $1$ at time $t_1$, and one type 0 individual (if there is one) becomes type 2 at time $t_2$.    
Fix a positive integer $L \leq N/2$.  For $i = 1,2$, let $Y_i(t)$ be the number of type $i$ individuals at time $t$ and let $B_i$ be the event that $L \leq \max_{t \geq 0} Y_i(t) \leq N/2$.  Then 
$$P(B_1 \cap B_2) \leq 2/L^2.$$
\end{Lemma}

\begin{Lemma}\label{qjbar}
Consider the model introduced in Proposition \ref{family}.  Assume $N \mu^{1 - 2^{-j}} \rightarrow \infty$ as $N \rightarrow \infty$.  We define the following three events:
\begin{enumerate}
\item Let $R_1$ be the event that eventually a type $j+1$ individual is born.

\item Let $R_2$ be the event that the maximum number of individuals of nonzero type at any time is between $\epsilon \mu^{-1 + 2^{-j}}$ and $N/2$.

\item Let $R_3$ be the event that all individuals still alive at time $\epsilon^{-1} \mu^{-1 + 2^{-j}}$ have type zero.
\end{enumerate}
Let ${\bar q}_{j+1} = P(R_1 \cap R_2 \cap R_3)$.  Then there exists a constant $C$, not depending on $\epsilon$, such that $q_{j+1} - C \epsilon \mu^{1-2^{-j}} \leq {\bar q}_{j+1} \leq q_{j+1}$.
\end{Lemma}

\begin{proof}
Because $q_{j+1} = P(R_1)$, the inequality ${\bar q}_{j+1} \leq q_{j+1}$ is immediate.  We need to show that $P(R_1 \cap (R_2^c \cup R_3^c)) \leq C \epsilon \mu^{1 - 2^{-j}}$.  Because $\epsilon^{-1} \mu^{-1 + 2^{-j}} \leq N$ for sufficiently large $N$, we have $P(R_3^c) \leq C \epsilon \mu^{1 - 2^{-j}}$ by (\ref{KolmPt}).  It remains to show that $P(R_1 \cap R_2^c) \leq C \epsilon \mu^{1 - 2^{-j}}$.

The probability that the number of individuals of nonzero type ever exceeds $N/2$ is at most $2/N \ll \epsilon \mu^{1 - 2^{-j}}$.  By (\ref{ELk}) and the fact that each type 1 individual experiences type 2 mutations at rate $\mu$, the expected number of type 2 mutations while there are $k$ individuals of nonzero type is at most $(k \mu)(1/k) = \mu$.  Therefore, the expected number of type 2 mutations while there are fewer than $\epsilon \mu^{-1+2^{-j}}$ individuals of nonzero type is at most $\epsilon \mu^{2^{-j}}$.  The probability that a given type 2 mutation has a type $j+1$ descendant is at most $C \mu^{1-2^{-(j-1)}}$ by Proposition \ref{family}.  It now follows, using Markov's Inequality, that the probability that some type 2 mutation that occurs while there are fewer than $\epsilon \mu^{-1+2^{-j}}$ individuals of nonzero type has a type $j+1$ descendant is at most $C \epsilon \mu^{2^{-j} + 1 - 2^{-(j-1)}} = C \epsilon \mu^{1 - 2^{-j}}$.  Thus, $P(R_1 \cap R_2^c) \leq C \epsilon \mu^{1 - 2^{-j}}$.  The result follows.
\end{proof}

We now define the events to which we will apply Lemma \ref{Poisson}.  Divide the interval $[0, T]$ into $M$ subintervals of equal length called $I_1, I_2, \dots, I_M$, where $M$ will tend to infinity with $N$.  Because type $m-j-1$ individuals experience type $m-j$ mutations at rate $\mu$, we can construct an inhomogeneous Poisson process $K$ on $[0, T]$ whose intensity at time $s$ is given by
\begin{equation}\label{intensity}
\frac{N \mu^{m-j} s^{m-j-1}}{(m-j-1)!} + \epsilon N \mu^{m-j} T^{m-j-1}
\end{equation}
such that on the event $G_N$, all the times of the type $m-j$ mutations before time $T$ are points of $K$.  Let $D_i$ be the event that there is a point of $K$ in the interval $I_i$.  Let $\xi_1, \xi_2, \dots, \xi_M$ be i.i.d. $\{0, 1\}$-valued random variables, independent of $K$ and the population process, such that $P(\xi_i = 1) = {\bar q}_{j+1}$ for all $i$, where ${\bar q}_{j+1}$ comes from Lemma \ref{qjbar}.  Let $A_i$ be the event that $D_i$ occurs, and one of the following occurs:
\begin{itemize}
\item The first point of $K$ in $I_i$ is the time of a type $m-j$ mutation, and the three events defined in Lemma \ref{qjbar} hold.  That is, the type $m-j$ mutation eventually has a type $m$ descendant, the maximum number of descendants that it has in the population at any future time is between $\epsilon \mu^{-1+2^{-j}}$ and $N/2$, and it has no descendants remaining a time $\epsilon^{-1} \mu^{-1 + 2^{-j}}$ after the mutation occurs.

\item There is no mutation at the time of the first point of $K$ in $I_i$, and $\xi_i = 1$.
\end{itemize}
Let $W = \sum_{i=1}^M {\bf 1}_{A_i}$ be the number of the events $A_i$ that occur, and let $\lambda = E[W]$.

\begin{Lemma}\label{Wlem}
We have $\limsup_{N \rightarrow \infty} |P(W = 0) - e^{-\lambda}| = 0.$
\end{Lemma}

\begin{proof}
Let $\beta_i$ be the set of all $j \leq M$ such that the distance between the intervals $I_i$ and $I_j$ is at most $\epsilon^{-1} \mu^{-1+2^{-j}}$.  Define $b_1$, $b_2$, and $b_3$ as in Lemma \ref{Poisson}.  We need to show that $b_1$, $b_2$, and $b_3$ all tend to zero as $N \rightarrow \infty$.

It is clear from properties of Poisson processes that the events $D_1, \dots, D_M$ are independent, and it is clear from the construction that $P(A_i|D_i) = {\bar q}_{j+1}$ for all $i$.  The events $A_1, \dots, A_M$ are not independent because mutations in two intervals $I_h$ and $I_i$ may have descendants alive at the same time.  However, if $I_i = [a, b]$, then the third event in Lemma \ref{qjbar} guarantees that whether or not $A_i$ has occurred is determined by time $b + \epsilon^{-1} \mu^{-1 + 2^{-j}}$, and therefore $A_i$ is independent of all $A_h$ with $h \notin \beta_i$.  It follows that $b_3 = 0$.

The length $|I_i|$ of the interval $I_i$ is $T/M$.  In view of (\ref{intensity}),
\begin{equation}\label{PDi}
P(D_i) \leq C N \mu^{m-j} T^{m-j-1} |I_i| = C N \mu^{m-j} T^{m-j}/M.
\end{equation}
Because (\ref{mujcond}) holds, we can apply Proposition \ref{family} to get ${\bar q}_{j+1} \leq q_{j+1} \leq C \mu^{1 - 2^{-j}}$.  Therefore, using also (\ref{latenote}),
$$P(A_i) = P(D_i) {\bar q}_{j+1} \leq \frac{C N \mu^{m-j+1-2^{-j}} T^{m-j}}{M} \leq \frac{C}{M}$$
for all $i$.  There are at most $2(1 + \epsilon^{-1} \mu^{-1 + 2^{-j}}/|I_i|) \leq C \epsilon^{-1} \mu^{-1 + 2^{-j}} M/T$ indices in $\beta_i$.  It follows that
\begin{align}\label{b1bound}
b_1 &\leq M \bigg( \frac{C \epsilon^{-1} \mu^{-1 + 2^{-j}} M}{T} \bigg) \bigg( \frac{C}{M} \bigg)^2 \nonumber \\
&\leq C \epsilon^{-1} \mu^{-1 + 2^{-j}} T^{-1} \nonumber \\
&\leq C \epsilon^{-1} \mu^{-1 + 2^{-j}} N^{1/(m-j)} \mu^{1 + (1-2^{-j})/(m-j)} \nonumber \\
&= C \epsilon^{-1} (N \mu^{1 + 2^{-j}(m-j-1)})^{1/(m-j)} \rightarrow 0
\end{align} 
as $N \rightarrow \infty$, using the second inequality in (\ref{mucond}).

To bound $b_2$, suppose $h \neq i$.  Suppose $D_h$ and $D_i$ both occur.  If the first points of the Poisson process in $I_h$ and $I_i$ are times of type $m-j$ mutations, then for $A_h \cap A_i$ to occur, the event $B_1 \cap B_2$ in Lemma \ref{2families} must occur with $L = \epsilon \mu^{-1 + 2^{-j}}$.  It follows that $$P(A_h \cap A_i|D_h \cap D_i) \leq \max\{2/(\epsilon \mu^{-1 + 2^{-j}})^2, {\bar q}_{j+1}^2 \} \leq C \epsilon^{-2} \mu^{2 - 2^{-(j-1)}}.$$  Therefore, using (\ref{PDi}), (\ref{latenote}), and the fact that $P(D_h \cap D_i) = P(D_h) P(D_i)$ by independence,
$$P(A_h \cap A_i) \leq P(D_h) P(D_i) P(A_h \cap A_i|D_h \cap D_i) \leq \bigg( \frac{C N \mu^{m-j} T^{m-j}}{M} \bigg)^2 (C \epsilon^{-2} \mu^{2 - 2^{-(j-1)}}) \leq \frac{C}{\epsilon^2 M^2}.$$  Thus, by reasoning as in (\ref{b1bound}), we get
$$b_2 \leq M \bigg( \frac{C \epsilon^{-1} \mu^{-1 + 2^{-j}} M}{T} \bigg) \bigg( \frac{C}{\epsilon^2 M^2} \bigg) \rightarrow 0$$ as $N \rightarrow \infty$, which completes the proof.
\end{proof}

\begin{Lemma}\label{sigmamlem}
Let $\sigma_m$ be the time of the first type $m-j$ mutation that will have a type $m$ descendant.  Then $$\lim_{N \rightarrow \infty} P(\sigma_m > T) = \exp \bigg( - \frac{t^{m-j}}{(m-j)!} \bigg).$$
\end{Lemma}

\begin{proof}
We claim there is a constant $C$, not depending on $\epsilon$, such that for sufficiently large $N$,
\begin{equation}\label{claim1}
\bigg| \lambda - \frac{t^{m-j}}{(m-j)!}\bigg| \leq C \epsilon,
\end{equation}
where $\lambda$ comes from Lemma \ref{Wlem}, and
\begin{equation} \label{claim2}
|P(W = 0) - P(\sigma_m > T)| \leq C \epsilon.
\end{equation}
The result follows from this claim by letting $\epsilon \rightarrow 0$ and applying Lemma \ref{Wlem}.

Recall that we have divided the interval $[0, T]$ into the subintervals $I_1, \dots, I_M$.  By letting $M$ tend to infinity sufficiently rapidly as $N$ tends to infinity, we can ensure that the expected number of points of the Poisson process $K$ that are in the same subinterval as some other point tends to zero as $N \rightarrow \infty$.  Therefore, $\sum_{i=1}^M P(D_i)$ is asymptotically equivalent to the expected number of points of $K$.  That is,
\begin{equation}\label{PDi2}
\sum_{i=1}^M P(D_i) \sim \int_0^T \frac{N \mu^{m-j} s^{m-j-1}}{(m-j-1)!} + \epsilon N \mu^{m-j} T^{m-j-1} \: ds = \frac{N \mu^{m-j} T^{m-j}}{(m-j)!} + \epsilon N \mu^{m-j} T^{m-j}.
\end{equation}
Now $$\lambda = \sum_{i=1}^M P(A_i) = {\bar q}_{j+1} \sum_{i=1}^M P(D_i),$$
so using Proposition \ref{family}, the second inequality in Lemma \ref{qjbar}, (\ref{PDi2}), and (\ref{latenote}),
\begin{equation}\label{lamup}
\limsup_{N \rightarrow \infty} \lambda \leq \limsup_{N \rightarrow \infty} \mu^{1 - 2^{-j}} \bigg( \frac{N \mu^{m-j} T^{m-j}}{(m-j)!} + \epsilon N \mu^{m-j} T^{m-j} \bigg) = \frac{t^{m-j}}{(m-j)!} + t^{m-j} \epsilon.
\end{equation}
Likewise, dropping the second term and using the first inequality in Lemma \ref{qjbar}, we get
\begin{equation}\label{lamlow}
\liminf_{N \rightarrow \infty} \lambda \geq \liminf_{N \rightarrow \infty} \: (1 - C \epsilon) \mu^{1 - 2^{-j}} \bigg( \frac{N \mu^{m-j} T^{m-j}}{(m-j)!} \bigg) = \frac{t^{m-j} (1 - C \epsilon)}{(m-j)!}.
\end{equation}
Equations (\ref{lamup}) and (\ref{lamlow}) imply (\ref{claim1}).

It remains to prove (\ref{claim2}).  The only way to have $W > 0$ and $\sigma_m > T$ is if for some $i$, there is a point of $K$ in $I_i$ that is not the time of a type $m-j$ mutation and $\xi_i = 1$.  On $G_N$, points of $K$ that are not mutation times occur at rate at most $2 \epsilon N \mu^{m-j} T^{m-j-1}$.  Because the Poisson process runs for time $T$ and $P(\xi_i = 1) = {\bar q}_{j+1} \leq C \mu^{1 - 2^{-j}}$ by Lemma \ref{qjbar} and Proposition \ref{family}, we have, using (\ref{latenote}),
\begin{equation}\label{way1}
P(W > 0 \mbox{ and }\sigma_m > T) \leq P(G_N^c) + C \epsilon N \mu^{m-j+1-2^{-j}} T^{m-j} \leq P(G_N^c) + C \epsilon.
\end{equation}

We can have $W = 0$ with $\sigma_m \leq T$ in two ways.  One possibility is that
two points of $K$ occur in the same subinterval, an event whose probability goes to zero if $M$ goes to infinity sufficiently rapidly with $N$.  The other possibility is that some type $m-j$ mutation before time $T$ could have a type $m$ descendant but fail to satisfy one of the other two conditions of Lemma \ref{qjbar}.  The probability of this event is at most
\begin{equation}\label{way2}
P(G_N^c) + C N \mu^{m-j} T^{m-j} (q_{j+1} - {\bar q}_{j+1}) \leq P(G_N^c) + C \epsilon N \mu^{m-j+1-2^{-j}} T^{m-j} \leq P(G_N^c) + C \epsilon
\end{equation}
by Lemma \ref{qjbar} and (\ref{latenote}).  Equation (\ref{claim2}) follows from (\ref{way1}), (\ref{way2}), and Lemma \ref{GNlem}.
\end{proof}

\begin{proof}[Proof of part 2 of Theorem \ref{bigmu}]
Recall the definition of $T$ from (\ref{defT}).
Define $\sigma_m$ to be the time of the first type $m-j$ mutation that will have a type $m$ descendant.  Then $\sigma_m \leq \tau_m$, and by Lemma \ref{sigmamlem}, it suffices to show that
\begin{equation}\label{ststs}
\lim_{N \rightarrow \infty} P(\sigma_m < T \mbox{ and }\tau_m - \sigma_m > \delta N^{-1/(m-j)} \mu^{-1 - (1-2^{-j})/(m-j)}) = 0
\end{equation}
for all $\delta > 0$.  The event in (\ref{ststs}) can only occur if some type $m-j$ mutation before time $T$ either fixates or takes longer than time $\delta N^{-1/(m-j)} \mu^{-1 - (1-2^{-j})/(m-j)}$ to disappear from the population.  By Lemma \ref{expk}, before time $T$ the expected rate of type $m-j$ mutations is at most $C N \mu^{m-j} T^{m-j-1}$, so the expected number of type $m-j$ mutations by time $T$ is at most $C N \mu^{m-j} T^{m-j}$.  Because the probability that a mutation fixates is $1/N$, the probability that some type $m-j$ mutation before time $T$ fixates is at most $C \mu^{m-j} T^{m-j}$, which goes to zero as $N \rightarrow \infty$ because $\mu T \rightarrow 0$ by (\ref{check1}).

Next, note that $\delta N^{-1/(m-j)} \mu^{-1 - (1-2^{-j})/(m-j)} \ll N$, which can be seen by dividing both sides by $N$ and observing that $\delta (N \mu)^{-1} (N \mu^{1-2^{-j}})^{-1/(m-j)} \rightarrow 0$ because $N \mu \rightarrow \infty$ and $N \mu^{1-2^{-j}} \rightarrow \infty$.  Therefore, for sufficiently large $N$, we can apply (\ref{KolmPt}) to show that the probability that a given mutation lasts longer than time $\delta N^{-1/(m-j)} \mu^{-1 - (1-2^{-j})/(m-j)}$ before disappearing or fixating is at most $C \delta^{-1} N^{1/(m-j)} \mu^{1 + (1-2^{-j})/(m-j)}$.  Thus, the probability that some mutation before time $T$ lasts this long is at most
\begin{align}
C \delta^{-1} N^{1/(m-j)} \mu^{1 + (1-2^{-j})/(m-j)} \cdot N \mu^{m-j} T^{m-j} &\leq C \delta^{-1} N^{1/(m-j)} \mu^{1 + (1-2^{-j})/(m-j)} \mu^{-(1 - 2^{-j})} t^{m-j} \nonumber \\
&= C \delta^{-1} (N \mu^{1 + (m-j-1)2^{-j}})^{1/(m-j)} t^{m-j} \rightarrow 0 \nonumber
\end{align}
by the second inequality in (\ref{mucond}), and (\ref{ststs}) follows.
\end{proof}

\section{Proof of part 3 of Theorem \ref{bigmu}}

Throughout this section, we assume
\begin{equation}\label{mup3}
\mu \sim A N^{-1/(1 + (m-j-1)2^{-j})}
\end{equation}
for some $j = 1, \dots, m-1$, as in part 3 of Theorem \ref{bigmu}.  Also, let $T = \mu^{-(1-2^{-j})} t$.  Then
\begin{equation}\label{muAeq}
\lim_{N \rightarrow \infty} N \mu^{m-j} T^{m-j} \mu^{1-2^{-j}} = \lim_{N \rightarrow \infty} N \mu^{1 + (m-j-1)2^{-j}} t^{m-j} = A^{1 + (m-j-1)2^{-j}} t^{m-j}.
\end{equation}
We first show that the number of individuals of type $m-j-1$ is approximately deterministic through time $T$.  

\begin{Lemma}\label{GNlem2}
Let $\epsilon > 0$.  Let $G_N(\epsilon)$ be the event that $$\max_{0 \leq s \leq T} \bigg|X_{m-j-1}(s) - \frac{N \mu^{m-j-1} s^{m-j-1}}{(m-j-1)!} \bigg| \leq \epsilon N \mu^{m-j-1} T^{m-j-1}.$$  Then $\lim_{N \rightarrow \infty} P(G_N(\epsilon)) = 1$.  \end{Lemma}

\begin{proof}
As in the proof of Lemma \ref{GNlem}, we need to check the conditions of Proposition \ref{detprop} with $m-j-1$ in place of $k$.  Because $\mu \rightarrow 0$ as $N \rightarrow \infty$, we have
\begin{equation}\label{muTzero}
\mu T = \mu^{2^{-j}} t \rightarrow 0
\end{equation}
as $N \rightarrow \infty$.  Also, using that $\mu \sim AN^{-1/(1 + (m-j-1)2^{-j})} \gg N^{-1/(1 + (m-j-2)2^{-j})}$, we have
$$N \mu^{m-j-1} T^{m-j-2} = N \mu^{m-j-1} \mu^{-(1-2^{-j})(m-j-2)} t^{m-j-2} = N \mu^{1 + (m-j-2)2^{-j}} t^{m-j-2} \rightarrow \infty$$ as $N \rightarrow \infty$.  Since $T \rightarrow \infty$ as $N \rightarrow \infty$, we also have $N \mu^{m-j-1} T^{m-j-1} \rightarrow \infty$ as $N \rightarrow \infty$, and the lemma follows.
\end{proof}

Although the number of type $m-j-1$ individuals is approximately deterministic, there are stochastic effects both from the number of type $m-j$ individuals in the population and from the time that elapses between the appearance of the type $m-j$ mutation that will have a type $m$ descendant and the birth of the type $m$ descendant.  Further complicating the proof is that because births and deaths occur at the same time in the Moran model, the fates of two type $m-j$ mutations that occur at different times are not independent, nor is the number of type $m-j$ individuals in the population independent of whether or not the type $m-j+1$ mutations succeed in producing a type $m$ descendant.  Our proof is very similar to the proof of Proposition 4.1 in \cite{DSS} and involves a comparison between the Moran model and a two-type branching process.  To carry out this comparison, we introduce five models.

\bigskip
\noindent {\bf Model 1}:  This will be the original model described in the introduction.

\bigskip
\noindent {\bf Model 2}:  This model is the same as Model 1 except that there are no type 1 mutations and no individuals of types $1, \dots, m-j-1$.  Instead, at times of an inhomogeneous Poisson process whose rate at time $s$ is $N \mu^{m-j} s^{m-j-1}/(m-j-1)! $, a type zero individual (if there is one) becomes type $m-j$.

\bigskip
\noindent {\bf Model 3}:  This model is the same as Model 2, except that type $m-j+1$ mutations are suppressed when there is another individual of type $m-j+1$ or higher already in the population.

\bigskip
\noindent {\bf Model 4}:  This model is the same as Model 3, except that two changes are made so that the evolution of type $m-j+1$ individuals and their offspring is decoupled from the evolution of the type $m-j$ individuals:
\begin{itemize}
\item Whenever there would be a transition that involves exchanging a type $m-j$ individual with an individual of type $k \geq m-j+1$, we instead exchange a randomly chosen type $0$ individual with a type $k$ individual.

\item At the times of type $m-j+1$ mutations, a randomly chosen type $0$ individual, rather than a type $m-j$ individual, becomes type $m-j+1$.
\end{itemize}

\smallskip
\noindent {\bf Model 5}:  This model is a two-type branching process with immigration.  Type $m-j$ individuals immigrate at times of an inhomogeneous Poisson process whose rate at time $s$ is $N \mu^{m-j} s^{m-j-1}/(m-j-1)!$.  Each individual gives birth at rate 1 and dies at rate 1, and type $m-j$ individuals become type $m$ at rate $\mu q_j$, where $q_j$ comes from Proposition \ref{family}.

\bigskip
For $i = 1, 2, 3, 4, 5$, let $Y_i(s)$ be the number of type $m-j$ individuals in Model $i$ at time $s$, and let $Z_i(s)$ be the number of individuals in Model $i$ at time $s$ of type $m-j+1$ or higher.  Let $r_i(s)$ be the probability that through time $s$, there has never been a type $m$ individual in Model $i$.  Note that $r_1(T) = P(\tau_m > T)$, so to prove part 3 of Theorem \ref{bigmu}, we need to calculate $\lim_{N \rightarrow \infty} r_1(T)$.  We will first find $\lim_{N \rightarrow \infty} r_5(T)$ and then bound $|r_i(T) - r_{i+1}(T)|$ for $i = 1, 2, 3, 4$.

\subsection{A two-type branching process with immigration}

Here we consider Model 5.  Our analysis is based on the following lemma concerning two-type branching processes, which is proved in section 2 of \cite{DSS}; see equation (2.4).

\begin{Lemma}\label{bplem}
Consider a continuous-time two-type branching process started with a single type $1$ individual.  Each type 1 individual gives birth and dies at rate one, and mutates to type 2 at rate $r$.  Let $f(t)$ be the probability that a type 2 individual is born by time $t$.  If $r$ and $t$ depend on $N$ with $r \rightarrow 0$ and $r^{1/2} t \rightarrow s$ as $N \rightarrow \infty$, then $$\lim_{N \rightarrow \infty} r^{-1/2} f(t) = \frac{1 - e^{-2s}}{1 + e^{-2s}}.$$
\end{Lemma}

\begin{Lemma}\label{r5lem}
We have
\begin{equation}\label{r5lemeq}
\lim_{N \rightarrow \infty} r_5(T) = \exp \bigg( - \frac{A^{1+(m-j-1)2^{-j}}}{(m-j-1)!} \int_0^t (t - s)^{m-j-1} \frac{1 - e^{-2s}}{1 + e^{-2s}} \: ds \bigg).
\end{equation}
\end{Lemma}

\begin{proof}
Let $g(w)$ be the probability that in Model 5, a type $m-j$ individual that immigrates at time $w$ has a type $m$ descendant by time $T$.  Because type $m-j$ individuals immigrate at times of an inhomogeneous Poisson process whose rate at time $w$ is $N \mu^{m-j} w^{m-j-1}/(m-j-1)!$, we have
\begin{equation}\label{r5eq}
r_5(T) = \exp \bigg( - \frac{1}{(m-j-1)!} \int_0^T N \mu^{m-j} w^{m-j-1} g(w) \: dw \bigg).
\end{equation}
Making the substitution $s = \mu^{1-2^{-j}} w$, we get
\begin{equation}\label{subs}
\int_0^T N \mu^{m-j} w^{m-j-1} g(w) \: dw = \int_0^t N \mu^{1 + (m-j-1)2^{-j}} s^{m-j-1} g(\mu^{-(1-2^{-j})} s) \mu^{-(1-2^{-j})} \: ds.
\end{equation}
As $N \rightarrow \infty$, we have $N \mu^{1 + (m-j-1)2^{-j}} \rightarrow A^{1 + (m-j-1)2^{-j}}$ by (\ref{mup3}).  Note also that $g(\mu^{-(1-2^{-j})}s) = f(\mu^{-(1-2^{-j})}(t-s))$, where $f$ is the function in Lemma \ref{bplem} when $r = \mu q_j$.  Also, by Proposition \ref{family}, $\mu q_j \sim \mu \cdot \mu^{1-2^{-(j-1)}} = (\mu^{1-2^{-j}})^2$, so $r^{-1/2} \sim \mu^{-(1-2^{-j})}$ and $r^{1/2} \mu^{-(1-2^{-j})} (t - s) \rightarrow t - s$ as $N \rightarrow \infty$.  Therefore, by Lemma \ref{bplem},
$$\lim_{N \rightarrow \infty} g(\mu^{-(1-2^{-j})}s) \mu^{-(1-2^{-j})} = \frac{1-e^{-2(t-s)}}{1+e^{-2(t-s)}}.$$
Using also (\ref{subs}) and the Dominated Convergence Theorem,
\begin{equation}\label{limint}
\lim_{N \rightarrow \infty} \int_0^T N \mu^{m-j} w^{m-j-1} g(w) \: dw = A^{1+ (m-j-1)2^{-j}} \int_0^t s^{m-j-1} \frac{1-e^{-2(t-s)}}{1+e^{-2(t-s)}} \: ds. 
\end{equation}
The result follows from (\ref{r5eq}) and (\ref{limint}) after interchanging the roles of $s$ and $t-s$.
\end{proof}

\subsection{Bounding the number of individuals of type $m-j$ and higher}

We begin with the following lemma, which bounds in all models the expected number of individuals in the models having type $m-j$ or higher. 

\begin{Lemma}\label{mutlem}
For $i = 1, 2, 3, 4, 5$, we have
\begin{equation}\label{maxexp}
\max_{0 \leq s \leq T} E[Y_i(s) + Z_i(s)] \leq C N \mu^{m-j} T^{m-j}.
\end{equation}
Also, for all five models, the expected number of type $m-j+1$ mutations by time $T$ is at most $C N \mu^{m-j+1} T^{m-j+1}$.
\end{Lemma}

\begin{proof}
Because each type $m-j$ individual experiences type $m-j+1$ mutations at rate $\mu$, the second statement of the lemma follows easily from the fact that $E[Y_i(s)] \leq C N \mu^{m-j} T^{m-j}$, which is a consequence of (\ref{maxexp}).

To prove (\ref{maxexp}), first note that because births and deaths occur at the same rate, in all five models $E[Y_i(s) + Z_i(s)]$ is the expected number of individuals of types $m-j$ and higher that appear up to time $s$ as a result of mutations, or immigration in the case of Model 5.  For $i = 2, 3, 5$, these mutation or immigration events occur at times of a rate $N \mu^{m-j} s^{m-j-1}/(m-j-1)!$ Poisson process (unless they are suppressed in Model 2 or 3 because no type zero individuals remain), so (\ref{maxexp}) holds.  In Model 1, the mutation rate depends on the number of type $m-j-1$ individuals, but (\ref{maxexp}) holds by Lemma \ref{expk}.

Model 4 is different because type $0$ rather than type $m-j$ individuals are replaced at the times of type $m-j+1$ mutations.  The above argument still gives $E[Y_4(s)] \leq C N \mu^{m-j} T^{m-j}$ for $s \leq T$ because type $m-j$ individuals give birth and die at the same rate.  Thus, the expected number of type $m-j+1$ mutations by time $T$ is at most $C N \mu^{m-j+1} T^{m-j+1}$.  It follows that $E[Z_4(s)] \leq C N \mu^{m-j+1} T^{m-j+1} \ll N \mu^{m-j} T^{m-j}$ for $s \leq T$, using the fact that $\mu T \rightarrow 0$ as $N \rightarrow \infty$ by (\ref{muTzero}).  Therefore, (\ref{maxexp}) holds for Model 4 as well.
\end{proof}

Lemma \ref{mutlem} easily implies the following bound on the maximum number of individuals of type $m-j$ or higher through time $T$.  The lemma below with $f(N) = 1/N$ implies that with probability tending to one as $N \rightarrow \infty$, the number of individuals of type $m-j$ or higher does not reach $N$ before time $T$.

\begin{Lemma}\label{maxbound}
Suppose $f$ is a function of $N$ such that $N \mu^{(m-j)2^{-j}} f(N) \rightarrow 0$ as $N \rightarrow \infty$.  Then for $i = 1, 2, 3, 4, 5$, as $N \rightarrow \infty$ we have
\begin{equation}\label{maxf}
\max_{0 \leq s \leq T} (Y_i(s) + Z_i(s)) f(N) \rightarrow_p 0.
\end{equation}
\end{Lemma}

\begin{proof}
Because individuals of type $m-j$ or higher give birth and die at the same rate, and they can appear but not disappear as a result of mutations, the process $(Y_i(s) + Z_i(s), 0 \leq s \leq T)$ is a nonnegative submartingale for $i = 1, 2, 3, 4, 5$.  By Doob's Maximal Inequality, for all $\delta > 0$,
\begin{equation}\label{delfN}
P \bigg( \max_{0 \leq s \leq T} (Y_i(s) + Z_i(s)) > \frac{\delta}{f(N)} \bigg) \leq \frac{E[Y_i(T) + Z_i(T)] f(N)}{\delta}.
\end{equation}
Since $N \mu^{m-j} T^{m-j} = N \mu^{(m-j) 2^{-j}} t^{m-j}$, equation (\ref{maxexp}) implies that if $N \mu^{(m-j)2^{-j}} f(N) \rightarrow 0$ as $N \rightarrow \infty$, then the right-hand side of (\ref{delfN}) goes to zero as $N \rightarrow \infty$ for all $\delta > 0$, which proves (\ref{maxf}).
\end{proof}

\subsection{Comparing Models 1 and 2}

In this subsection, we establish the following result which controls the difference between Model 1 and Model 2.  The advantage to working with Model 2 rather than Model 1 is that the randomness in the rate of the type $m-j$ mutations is eliminated.

\begin{Lemma}\label{r1r2}
We have $\lim_{N \rightarrow \infty} |r_1(T) - r_2(T)| = 0$.
\end{Lemma}

\begin{proof}
Lemma \ref{maxbound} with $f(N) = 1/N$ implies that with probability tending to one as $N \rightarrow \infty$, up to time $T$ there is always at least one type $0$ individual in Model 2, so hereafter we will make this assumption.  In this case, a type $m-j$ individual replaces a randomly chosen type 0 individual in Model 2 at times of a Poisson process $K$ whose rate at time $s$ is $N \mu^{m-j} s^{m-j-1}/(m-j-1)!$.  We will first compare Model 2 to another model called Model $2'$, which will be the same as Model 2 except that type $m-j$ individuals arrive at times of a Poisson process $K'$ whose rate at time $s$ is $\max\{0, N \mu^{m-j} s^{m-j-1}/(m-j-1)! - \epsilon N \mu^{m-j} T^{m-j-1} \}$, where $\epsilon > 0$ is fixed.

Models 2 and $2'$ can be coupled so that births and deaths occur at the same times in both models, and each point of $K'$ is also a point of $K$.  Consequently, a coupling can be achieved so that if an individual has type $k \geq m-j$ in Model $2'$, then it also has type $k$ in Model $2$.  With such a coupling, the only individuals whose types are different in the two models are those descended from 
individuals that in Model 2 became type $m-j$ at a time that is in $K$ but not $K'$.  The rate of points in $K$ but not $K'$ is bounded by $\epsilon N \mu^{m-j}  T^{m-j-1}$.  The probability that a given type $m-j$ individual has a type $m$ descendant is at most $C \mu^{1-2^{-j}}$ by Proposition \ref{family}.  Therefore, the probability that there is a type $m$ individual in Model 2 but not Model $2'$ before time $T$ is bounded by
\begin{equation}\label{epgap}
\epsilon N \mu^{m-j} T^{m-j} \cdot C \mu^{1-2^{-j}} \leq C \epsilon,
\end{equation}
using (\ref{muAeq}).  Therefore, letting $r_{2'}(T)$ denote the probability that there is no type $m$ individual in Model $2'$ by time $T$,
\begin{equation} \label{r22'}
|r_2(T) - r_{2'}(T)| \leq C \epsilon.
\end{equation}

We now compare Model 1 and Model $2'$.  These models can be coupled so that births and deaths in the two models happen at the same times and, on $G_N(\epsilon)$, there is a type $m-j$ mutation in Model 1 at all of the times in $K'$.  This coupling can therefore achieve the property that on $G_N(\epsilon)$, any individual of type $k \geq m-j$ in Model $2'$ also has type $k$ in Model 1.  The only individuals in Model 1 of type $k \geq m-j$ that do not have the same type in Model $2'$ are those descended from individuals that became type $m-j$ at a time that is not in $K'$.  On $G_N(\epsilon)$, the rate of type $m-j$ mutations at times not in $K'$ is bounded by $2 \epsilon N \mu^{m-j} T^{m-j-1}$.  Therefore, by the same calculation made in (\ref{epgap}), the probability that $G_N(\epsilon)$ occurs and that Model 1 but not Model $2'$ has a type $m$ descendant by time $T$ is at most $C \epsilon$.  This bound and Lemma \ref{GNlem2} give
\begin{equation} \label{r12'}
|r_1(T) - r_{2'}(T)| \leq C \epsilon.
\end{equation}
The result follows from (\ref{r22'}) and (\ref{r12'}) after letting $\epsilon \rightarrow 0$.
\end{proof}

\subsection{Comparing Models 2 and 3}

In this subsection, we establish the following lemma.

\begin{Lemma}\label{r2r3}
We have $\lim_{N \rightarrow \infty} |r_2(T) - r_3(T)| = 0$.
\end{Lemma}

The advantage to working with Model 3 rather than Model 2 is that in Model 3, descendants of only one type $m-j+1$ mutation can be present in the population at a time.  As a result, each type $m-j+1$ mutation independently has probability $q_j$ of producing a type $m$ descendant.  With Model 2, there could be dependence between the outcomes of different type $m-j+1$ mutations whose descendants overlap in time.

The only difference between Model 2 and Model 3 is that some type $m-j+1$ mutations are suppressed in Model 3.  Therefore, it is easy to couple Model 2 and Model 3 so that until there are no type 0 individuals remaining in Model 2, the type of the $i$th individual in Model 2 is always at least as large as the type of the $i$th individual in Model 3, with the only discrepancies involving individuals descended from a type $m-j+1$ mutation that was suppressed in Model 3.  Because Lemma \ref{maxbound} with $f(N) = 1/N$ implies that the probability that all type zero individuals disappear by time $T$ goes to zero as $N \rightarrow \infty$, Lemma \ref{r2r3} follows from the following result.

\begin{Lemma}
In Model 2, the probability that some type $m-j+1$ mutation that occurs while there is another individual of type $m-j+1$ or higher in the population has a type $m$ descendant tends to zero as $N \rightarrow \infty$.
\label{supp24}
\end{Lemma}

\begin{proof}
By Lemma \ref{mutlem}, the expected number of type $m-j+1$ mutations by time $T$ is at most $C N \mu^{m-j+1} T^{m-j+1}$.  By (\ref{ET}), the expected amount of time, before time $T$, that there is an individual in the population of type $m-j+1$ or higher is at most $C N \mu^{m-j+1} T^{m-j+1} (\log N)$.

By Lemma \ref{maxbound} with $f(N) = 1/(N \mu^{m-j} T^{m-j} \log N)$, the probability that the number of type $m-j$ individuals stays below $N \mu^{m-j} T^{m-j} \log N$ until time $T$ tends to one as $N \rightarrow \infty$.  On this event, the expected number of type $m-j+1$ mutations by time $T$ while there is another individual in the population of type $m-j+1$ or higher is at most $$h_N = (C N \mu^{m-j+1} T^{m-j+1} \log N) (N \mu^{m-j} T^{m-j} \log N) \mu.$$  The probability that a given such mutation produces a type $m$ descendant is $q_j \leq C \mu^{1-2^{-(j-1)}}$ by Proposition \ref{family}, so the probability that at least one such mutation produces a type $m$ descendant is at most $$h_N q_j \leq C (\mu T (\log N)^2) [N \mu^{m-j} T^{m-j} \mu^{1-2^{-j}}]^2.$$  Because $\mu T (\log N)^2 = \mu^{2^{-j}} (\log N)^2 \rightarrow 0$ as $N \rightarrow \infty$ and $N \mu^{m-j} T^{m-j} \mu^{1-2^{-j}}$ stays bounded as $N \rightarrow \infty$ by (\ref{muAeq}), the lemma follows.
\end{proof}

\subsection{Comparing Models 3 and 4}

In both Model 3 and Model 4, each type $m-j+1$ mutation independently has probability $q_j$ of producing a type $m$ descendant.  The advantage to Model 4 is that whether or not a given type $m-j+1$ mutation produces a type $m$ descendant is decoupled from the evolution of the number of type $m-j$ individuals.

We first define a more precise coupling between Model 3 and Model 4.  We will assume throughout the construction that there are fewer than $N/2$ individuals in each model with type $m-j$ or higher.  Eventually this assumption will fail, but by Lemma \ref{maxbound}, the assumption is valid through time $T$ with probability tending to one as $N \rightarrow \infty$, which is sufficient for our purposes.

For both models, the $N$ individuals will be assigned labels $1, \dots, N$ in addition to their types.  Let $L$ be a Poisson process of rate $N$ on $[0, \infty)$, and let $I_1, I_2, \dots$ and $J_1, J_2, \dots$ be independent random variables, uniformly distributed on $\{1, \dots, N\}$.  Let $K$ be an inhomogeneous Poisson process on $[0, \infty)$ whose rate at time $s$ is $N \mu^{m-j} s^{m-j-1}/(m-j-1)!$, and let $L_1, \dots, L_N$ be independent rate $\mu$ Poisson processes on $[0, \infty)$.  In both models, if $s$ is a point of $K$, then at time $s$ we choose an individual at random from those that have type $0$ in both models to become type $m-j$.  Birth and death events occur at the times of $L$.  At the time of the $m$th point of $L$, in both models we change the type of the individual labeled $I_m$ to the type of the individual labeled $J_m$.  In Model 4, if $I_m$ has type $m-j$ and $J_m$ has type $k \geq m-j+1$, then we choose a type 0 individual to become type $m-j$ to keep the number of type $m-j$ individuals constant.  Likewise, in Model 4, if $I_m$ has type $k \geq m-j+1$ and $J_m$ has type $m-j$, then we choose a type $m-j$ individual to become type 0.  In both models, the individual labeled $i$ experiences mutations at times of $L_i$, with the exceptions that type 0 individuals never get mutations and mutations of type $m-j$ individuals are suppressed when there is already an individual of type $m-j+1$ or higher in the population.  Also, in Model 4, if $s$ is a point of $L_i$ and the individual labeled $i$ has type $m-j$ at time $s-$, then in addition to changing the type of the individual labeled $i$, we choose a type 0 individual to become type $m-j$ so that the number of type $m-j$ individuals stays constant.

Note that by relabeling the individuals, if necessary, after each transition, we can ensure that for all $s \geq 0$, at time $s$ there are $\min\{Y_3(s), Y_4(s)\}$ integers $i$ such that the individual labeled $i$ has type $m-j$ in both models.  The rearranging can be done so that no individual has type $m-j$ in one of the models and type $m-j+1$ or higher in the other.  Also, with this coupling, if a type $m-j+1$ mutation occurs at the same time in both models, descendants of this mutation will have the same type in both models.  In particular, if the mutation has a type $m$ descendant in one model, it will have a type $m$ descendant in the other.

Let $W(s) = Y_3(s) - Y_4(s)$, which is the difference between the number of type $m-j$ individuals in Model 3 and the number of type $m-j$ individuals in Model 4.  There are three types of events that can cause the process $(W(s), 0 \leq s \leq T)$ to jump:
\begin{itemize}
\item When a type $m-j$ individual experiences a mutation in Model 3 and becomes type $m-j+1$, there is no change to the number of type $m-j$ individuals in Model 4.  At time $s$, such changes occur at rate either 0 or $\mu Y_3(s)$, depending on whether or not there is already an individual in Model 3 of type $m-j+1$ or higher.

\item When one of the individuals that is type $m-j$ in one process but not the other experiences a birth or death, the $W$ process can increase or decrease by one.  If $Y_3(s) > Y_4(s)$, then at time $s$, both increases and decreases are happening at rate $|W(s)|(N - |W(s)|)/N$ because the $W$ process changes unless the other individual involved in the exchange also has type $m-j$ in Model 3 but not Model 4.  If $Y_4(s) > Y_3(s)$, then increases and decreases are each happening at rate $|W(s)|(N - |W(s)| - Z_4(s))/N$ because in Model 4, transitions exchanging a type $m-j$ individual with an individual of type $m-j+1$ or higher are not permitted.

\item The number of type $m-j$ individuals changes in Model 3 but not Model 4 when there is an exchange involving one of the individuals that has type $m-j$ in both models and one of the individuals that has type $m-j+1$ or higher in Model 4.  Changes in each direction happen at rate $Z_4(s) \min\{Y_3(s), Y_4(s)\}/N$.
\end{itemize}
Therefore, the process $(W(s), 0 \leq s \leq T)$ at time $s$ is increasing by one at rate $\lambda(s)$ and decreasing by one at rate $\lambda(s) + \gamma(s)$, where
\begin{equation}\label{gamdef}
0 \leq \gamma(s) \leq \mu Y_3(s)
\end{equation}
and
\begin{equation}\label{lamdef}
\lambda(s) = \frac{|W(s)|(N - |W(s)| - Z_4(s) {\bf 1}_{\{Y_4(s) > Y_3(s)\}})}{N} + \frac{Z_4(s) \min\{Y_3(s), Y_4(s)\}}{N}.
\end{equation}
The next lemma bounds the process $(W(s), 0 \leq s \leq T)$.

\begin{Lemma}\label{EKlem}
For $0 \leq s \leq t$, let $$W_N(s) = \frac{1}{N \mu^{(m-j)2^{-j}}} W(s \mu^{-(1-2^{-j})}).$$  Then as $N \rightarrow \infty$,
\begin{equation}\label{maxWN}
\max_{0 \leq s \leq t} |W_N(s)| \rightarrow_p 0.
\end{equation}
\end{Lemma}

\begin{proof}
The proof is similar to the proof of Lemma 4.6 in \cite{DSS}.  We use Theorem 4.1 in chapter 7 of \cite{ek86} to show that the processes $(W_N(s), 0 \leq s \leq t)$ converge as $N \rightarrow \infty$ to a diffusion $(X(s), 0 \leq s \leq t)$ which satisfies the stochastic differential equation
\begin{equation}\label{SDE}
d X(s) = b(X(s)) + a(X(s)) \: dB(s)
\end{equation}
with $b(x) = 0$ and $a(x) = 2A^{- 1 - (m-j-1)2^{-j}} |x|$ for all $x$, where $A$ is the constant from (\ref{mup3}).  The Yamada-Watanabe Theorem (see, for example, (3.3) on p. 193 of \cite{durrsc}) gives pathwise uniqueness for this SDE, which implies that the associated martingale problem is well-posed.

For all $N$ and all $s \in [0, t]$, define
$$B_N(s) = - \frac{1}{N \mu^{(m-j)2^{-j}}} \int_0^s \frac{\gamma(r \mu^{-(1-2^{-j})})}{\mu^{1-2^{-j}}} \: dr = - \frac{1}{N \mu^{1 + (m-j-1)2^{-j}}} \int_0^s \gamma(r \mu^{-(1-2^{-j})}) \: dr$$ and
$$A_N(s) = \frac{1}{(N \mu^{(m-j)2^{-j}})^2 \mu^{1-2^{-j}}} \int_0^s \big(2 \lambda(r \mu^{-(1-2^{-j})}) + \gamma(r \mu^{-(1-2^{-j})})\big) \: dr.$$  At time $s$, the process $(W_N(s), 0 \leq s \leq t)$ experiences positive jumps by $1/(N \mu^{(m-j)2^{-j}})$ at rate $\lambda(s \mu^{-(1-2^{-j})}) \mu^{-(1-2^{-j})}$ and negative jumps by the same amount at the slightly larger rate $(\lambda(s \mu^{-(1-2^{-j})}) + \gamma(s \mu^{-(1-2^{-j})})) \mu^{-(1-2^{-j})}$.  Therefore, letting $M_N(s) = W_N(s) - B_N(s)$, the processes $(M_N(s), 0 \leq s \leq t)$ and $(M_N^2(s) - A_N(s), 0 \leq s \leq t)$ are martingales.  We claim that as $N \rightarrow \infty$, 
\begin{equation}\label{infmean}
\sup_{0 \leq s \leq t} |B_N(s)| \rightarrow_p 0
\end{equation}
and
\begin{equation}\label{infvar}
\sup_{0 \leq s \leq t} \bigg| A_N(s) - 2A^{-1 - (m-j-1)2^{-j}} \int_0^s |W_N(r)| \: dr \bigg| \rightarrow_p 0.
\end{equation}
The results (\ref{infmean}) and (\ref{infvar}) about the infinitesimal mean and variance respectively enable us to deduce from Theorem 4.1 in chapter 7 of \cite{ek86} that as $N \rightarrow \infty$, the processes $(W_N(s), 0 \leq s \leq T)$ converge in the Skorohod topology to a process $(X(s), 0 \leq s \leq T)$ satisfying (\ref{SDE}).  Because $W_N(0) = 0$ for all $N$, we have $X(0) = 0$, and therefore $X(s) = 0$ for $0 \leq s \leq T$.  The result (\ref{maxWN}) follows.

To complete the proof, we need to establish (\ref{infmean}) and (\ref{infvar}).  Equation (\ref{gamdef}) and Lemma \ref{maxbound} with $f(N) = t/(N \mu^{(m-j-1)2^{-j}})$ imply that as $N \rightarrow \infty$, $$\sup_{0 \leq s \leq t} |B_N(s)| \leq \frac{t}{N \mu^{(m-j-1)2^{-j}}} \max_{0 \leq s \leq T} Y_3(s) \rightarrow_p 0,$$ which proves (\ref{infmean}).

To prove (\ref{infvar}), note that
\begin{align}
A_N(s) - 2&A^{-1-(m-j-1)2^{-j}} \int_0^s |W_N(r)| \: dr \nonumber \\
&= \int_0^s \frac{ 2 \lambda(r \mu^{-(1-2^{-j})}) + \gamma(r \mu^{-(1-2^{-j})})}{(N \mu^{(m-j)2^{-j}})^2 \mu^{1-2^{-j}}} - \frac{2A^{-1-(m-j-1)2^{-j}}|W(r \mu^{-1-2^{-j}})|}{N \mu^{(m-j)2^{-j}}} \: dr. \nonumber
\end{align}
It therefore follows from (\ref{gamdef}) and (\ref{lamdef}) that
\begin{align}\label{4terms}
\sup_{0 \leq s \leq t} &\bigg| A_N(s) - 2A^{-1 - (m-j-1)2^{-j}} \int_0^s |W_N(r)| \: dr \bigg| \nonumber \\
&\leq \sup_{0 \leq s \leq t} \int_0^s \bigg| \frac{2}{(N \mu^{(m-j)2^{-j}})^2 \mu^{1-2^{-j}}} - \frac{2 A^{-1 - (m-j-1)2^{-j}}}{N \mu^{(m-j)2^{-j}}} \bigg| |W(r \mu^{-(1-2^{-j})})| \: dr \nonumber \\
&\: \: + \sup_{0 \leq s \leq t} \int_0^s \frac{2 W(r \mu^{-(1-2^{-j})})^2 + 2 |W(r \mu^{-(1-2^{-j})})| Z_4(r \mu^{-(1-2^{-j})})}{N (N \mu^{(m-j)2^{-j}})^2 \mu^{1-2^{-j}}} \: dr \nonumber \\
&\: \: + \sup_{0 \leq s \leq t} \int_0^s \frac{2 Z_4(r \mu^{-(1-2^{-j})}) \min\{Y_3(r \mu^{-(1-2^{-j})}), Y_4(r \mu^{-(1-2^{-j})})\}}{N (N \mu^{(m-j)2^{-j}})^2 \mu^{1-2^{-j}}} \: dr \nonumber \\
&\: \: + \sup_{0 \leq s \leq t} \int_0^s \frac{\mu Y_3(r \mu^{-(1 - 2^{-j})})}{(N \mu^{(m-j)2^{-j}})^2 \mu^{1-2^{-j}}} \: dr.
\end{align}
We need to show that the four terms on the right-hand side of (\ref{4terms}) each converge in probability to zero.  Because $t$ is fixed, in each case it suffices to show that the supremum of the integrand over $r \in [0, t]$ converges in probability to zero as $N \rightarrow \infty$.  We have
\begin{align}
\sup_{0 \leq s \leq T} &\bigg| \frac{2}{(N \mu^{(m-j)2^{-j}})^2 \mu^{1-2^{-j}}} - \frac{2A^{-1-(m-j-1)2^{-j}}}{N \mu^{(m-j)2^{-j}}} \bigg| \: |W(s)| \nonumber \\
&= \sup_{0 \leq s \leq T} \bigg| \frac{2}{N \mu^{1 + (m-j-1)2^{-j}}} - \frac{2}{A^{1+(m-j-1)2^{-j}}} \bigg| \cdot \frac{|W(s)|}{N \mu^{(m-j)2^{-j}}} \rightarrow_p 0 \nonumber
\end{align}
by Lemma \ref{maxbound} because $|W(s)| \leq \max\{Y_3(s), Y_4(s)\}$ and the first factor goes to zero as $N \rightarrow \infty$ by (\ref{mup3}).  Thus, the first term in (\ref{4terms}) converges in probability to zero.  Also, $N \mu^{1-2^{-j}} \rightarrow \infty$ as $N \rightarrow \infty$, so Lemma \ref{maxbound} gives
\begin{align}
\sup_{0 \leq s \leq T} &\frac{W(s)^2 + |W(s)| Z_4(s)}{N (N \mu^{(m-j)2^{-j}})^2 \mu^{1-2^{-j}}} \nonumber \\
&= \sup_{0 \leq s \leq T} \bigg( \frac{|W(s)|}{N \mu^{(m-j)2^{-j}} (N \mu^{1-2^{-j}})^{1/2}} \bigg) \bigg( \frac{|W(s)| + Z_4(s)}{N \mu^{(m-j)2^{-j}} (N \mu^{1 - 2^{-j}})^{1/2}} \bigg) \rightarrow_p 0, \nonumber
\end{align}
which is enough to control the second term in (\ref{4terms}).  The same argument works for the third term, using $Z_4(s) Y_4(s)$ in the numerator of the left-hand side in place of $W(s)^2 + |W(s)| Z_4(s)$.  Finally,
$$\sup_{0 \leq s \leq T} \frac{\mu Y_3(s)}{(N \mu^{(m-j)2^{-j}})^2 \mu^{1-2^{-j}}} = \frac{\mu Y_3(s)}{N \mu^{(m-j)2^{-j}}} \cdot \frac{1}{N \mu^{1 + (m-j-1)2^{-j}}} \rightarrow_p 0$$
by Lemma \ref{maxbound} because $\mu \rightarrow 0$ as $N \rightarrow \infty$ and $N \mu^{1 + (m-j-1)2^{-j}}$ is bounded away from zero as $N \rightarrow \infty$ by (\ref{mup3}).  Therefore, the fourth term on the right-hand side of (\ref{4terms}) converges in probability to zero, which completes the proof of (\ref{infvar}).
\end{proof}

\begin{Lemma}\label{endlem}
In both Model 3 and Model 4, the probability that there is a type $m-j+1$ mutation before time $T$ that has a type $m$ descendant born after time $T$ converges to zero as $N \rightarrow \infty$.
\end{Lemma}

\begin{proof}
The same argument works for both models.  Let $\epsilon > 0$.  By Lemma \ref{mutlem}, the expected number of type $m-j+1$ mutations by time $T$ is at most $N \mu^{m-j+1} T^{m-j+1}$.  Since $N \mu^{1-2^{-j}} \rightarrow \infty$ as $N \rightarrow \infty$, we have $\epsilon T \ll N$.  Therefore, by (\ref{KolmPt}), the probability that a given mutation stays in the population for a time at least $\epsilon T$ before dying out or fixating is at most $C/(\epsilon T)$.  It follows that the probability that some type $m-j+1$ mutation before time $T$ lasts for a time at least $\epsilon T$ is at most
$$C \epsilon^{-1} N \mu^{m-j+1} T^{m-j} \leq C \epsilon^{-1} N \mu^{1 + (m-j)2^{-j}} \rightarrow 0$$
as $N \rightarrow \infty$ by (\ref{mup3}).  Thus, with probability tending to one as $N \rightarrow \infty$, all type $m-j+1$ mutations that have a descendant alive at time $T$ originated after time $(1 - \epsilon)T$.  

Arguing as above, the expected number of type $m-j+1$ mutations between times $(1 - \epsilon)T$ and $T$ is at most $\epsilon N \mu^{m-j+1} T^{m-j+1}$, and the probability that a given such mutation has a type $m$ descendant is $q_j \leq C \mu^{1-2^{-(j-1)}}$ by Proposition \ref{family}.  Thus, the probability that some type $m-j+1$ mutation between times $(1 - \epsilon)T$ and $T$ has a type $m$ descendant is at most
\begin{equation}\label{CNmu}
C \epsilon N \mu^{m-j+1} T^{m-j+1} \mu^{1-2^{-(j-1)}} \leq C \epsilon N \mu^{1 + (m-j-1)2^{-j}} \leq C \epsilon
\end{equation}
by (\ref{mup3}).  The lemma follows by letting $\epsilon \rightarrow 0$.
\end{proof}

\begin{Lemma}\label{r3r4}
We have $\lim_{N \rightarrow \infty} |r_3(T) - r_4(T)| = 0$.
\end{Lemma}

\begin{proof}
For $i = 3, 4$, let $D_i$ be the event that no type $m-j+1$ mutation that occurs before time $T$ has a type $m$ descendant.  By Lemma \ref{endlem}, it suffices to show that
\begin{equation}\label{D3D4}
\lim_{N \rightarrow \infty} |P(D_3) - P(D_4)| = 0.
\end{equation}

Recall that Model 3 and Model 4 are coupled so that when a type $m-j+1$ mutation occurs at the same time in both models, it will have a type $m$ descendant in one model if and only if it has a type $m$ descendant in the other.  Therefore, $|P(D_3) - P(D_4)|$ is at most the probability that some type $m-j+1$ mutation that occurs in one process but not the other has a type $m$ descendant.  There are two sources of type $m-j+1$ mutations that occur in one process but not the other.  Some type $m-j+1$ mutations are suppressed in one model but not the other because there is already an individual of type $m-j+1$ or higher in the population.  That the probability of some such mutation having a type $m$ descendant goes to zero follows from the argument used to prove Lemma \ref{supp24}, which is also valid for Model 3 and Model 4.  The other type $m-j+1$ mutations that appear in one process but not the other occur when one of the $|W(s)|$ individuals that has type $m-j$ in one model but not the other gets a mutation.  Let $\epsilon > 0$.  By Lemma \ref{EKlem}, for sufficiently large $N$, $$P \bigg( \max_{0 \leq s \leq T} |W(s)| \leq \epsilon N \mu^{(m-j) 2^{-j}} \bigg) > 1 - \epsilon.$$  Therefore, on an event of probability at least $1 - \epsilon$, the expected number of type $m-j+1$ mutations that occur in one model but not the other and have a type $m$ descendant is at most $$\epsilon N \mu^{(m-j) 2^{-j}} q_j \leq C \epsilon N \mu^{1 + (m-j-1)2^{-j}} \leq C \epsilon$$ by Proposition \ref{family} and (\ref{mup3}).  The result follows by letting $\epsilon \rightarrow 0$.
\end{proof}

\subsection{Comparing Models 4 and 5}

In both Model 4 and Model 5, type $m-j$ individuals appear at times of a Poisson process whose rate at time $s$ is $N \mu^{m-j} s^{m-j-1}/(m-j-1)!$.  In both models, type $m-j$ individuals experience mutations that will lead to type $m$ descendants at rate $\mu q_j$.  The two models differ in the following three ways:
\begin{itemize}
\item In Model 4, some type $m-j+1$ mutations are suppressed because there is another individual of type $m-j+1$ or higher already in the population.

\item In Model 4, some time elapses between the time of the type $m-j+1$ mutation that will produce a type $m$ descendant, and the time that the type $m-j+1$ descendant appears.

\item In Model 4, when there are $k$ individuals of type $m-j$ and $\ell$ individuals of type $m-j+1$ or higher, the rate at which the number of type $m-j$ individuals increases (or decreases) by one is $k (N-\ell)/N$ because the number of type $m-j$ individuals changes only when a type $m-j$ individual is exchanged with a type 0 individual.  This rate is simply $k$ in Model 5.  An additional complication is that the factor $(N - \ell)/N$ is not independent of whether previous type $m-j+1$ mutations are successful in producing type $m$ descendants.
\end{itemize}

We prove Lemma \ref{r4r5} below by making three modifications to Model 4 to eliminate these differences, and then comparing the modified model to Model 5.  Lemmas \ref{r5lem}, \ref{r1r2}, \ref{r2r3}, \ref{r3r4}, and \ref{r4r5} immediately imply part 3 of Proposition \ref{bigmu}.

\begin{Lemma}\label{r4r5}
We have $\lim_{N \rightarrow \infty} |r_4(T) - r_5(T)| = 0$.
\end{Lemma}

\begin{proof}
We obtain Model $4'$ from Model 4 by making the following modifications.  First, whenever a type $m-j+1$ mutation is suppressed in Model 4 because there is another individual in the population of type $m-j+1$ or higher, in Model $4'$ we add a type $m$ individual with probability $q_j$.  Second, whenever a type $m-j+1$ mutation occurs in Model 4 that will eventually produce a type $m$ descendant, we change the type of the mutated individual in Model $4'$ to type $m$ immediately.  Third, for every type $m-j+1$ mutation in Model $4'$, including the events that produce a type $m$ individual that were added in the first modification, if there are $\ell$ individuals of type $m-j$ or higher in the population, then we suppress the mutation with probability $\ell/N$.  This means that at all times, every type $m-j$ individual in Model $4'$ experiences a mutation that will produce a type $m$ descendant at rate $\mu q_j (N-\ell)/N$, while new type $m-j$ individuals appear and disappear at rate $k(N-\ell)/N$.  Note that the number of type $m-j$ individuals is always the same in Model $4'$ as in Model 4.  Let $r_{4'}(T)$ be the probability that there is a type $m$ individual in Model $4'$ by time $T$.

Lemma \ref{supp24}, whose proof is also valid for Model $4'$, implies that with probability tending to one as $N \rightarrow \infty$, the first modification above does not cause a type $m$ individual to be added to Model $4'$ before time $T$.  Lemma \ref{endlem} implies this same result for the second modification.  As for the third modification, let $\epsilon > 0$, and let $D_N$ be the event that the number of individuals of type $m-j$ or higher in Model 4 stays below $\epsilon N$ through time $T$.  By Lemma \ref{maxbound}, we have $\lim_{N \rightarrow \infty} P(D_N) = 1$.  By Lemma \ref{mutlem}, the expected number of type $m-j+1$ mutations by time $T$ is at most $C N \mu^{m-j+1} T^{m-j+1}$.  On $D_N$, we always have $\ell/N < \epsilon$, so the probability that $D_N$ occurs and a type $m-j+1$ mutation that produces a type $m$ descendant in Model 4 gets suppressed in Model $4'$ is at most $C N \mu^{m-j+1} T^{m-j+1} \cdot q_j \epsilon \leq C \epsilon$, using (\ref{CNmu}) and Proposition \ref{family}.  Thus,
\begin{equation}\label{r44'}
\limsup_{N \rightarrow \infty} |r_4(T) - r_{4'}(T)| < \epsilon.
\end{equation}

It remains to compare Model $4'$ and Model 5.  In Model 5, when there are $k$ type $m-j$ individuals, the rates that type $m-j$ individuals appear, disappear, and give rise to a type $m$ individual are $k$, $k$, and $k \mu q_j$ respectively, as compared with $k(N-\ell)/N$, $k(N-\ell)/N$, and $k \mu q_j (N-\ell)/N$ respectively in Model $4'$.  Consequently, Model $4'$ is equivalent to Model 5 slowed down by a factor of $(N-\ell)/N$, which on $D_N$ stays between $1-\epsilon$ and $1$.  We can obtain a lower bound for $r_{4'}(T)$ by considering Model 5 run all the way to time $T$, so $r_{4'}(T) \geq r_5(T)$.  An upper bound for $r_{4'}(T)$ on $D_N$ is obtained by considering Model 5 run only to time $T(1 - \epsilon)$, so $r_4'(T) \leq r_5((1 - \epsilon)T) + P(D_N^c)$.  Now $\lim_{N \rightarrow \infty} r_5((1 - \epsilon)T)$ is given by the right-hand side of (\ref{r5lemeq}) with $(1-\epsilon)t$ in place of $t$.  Therefore, by letting $N \rightarrow \infty$ and then $\epsilon \rightarrow 0$, we get
$$\lim_{N \rightarrow \infty} |r_{4'}(T) - r_5(T)| = 0,$$ which, combined with (\ref{r44'}), proves the lemma.
\end{proof}

\bigskip
\noindent {\bf \Large Acknowledgments}
\bigskip

\noindent The author thanks Rick Durrett for many helpful discussions regarding this work.  He also thanks Rinaldo Schinazi for a discussion related to section \ref{powersec}, and a referee for comments about the presentation of the paper.


\begin{thebibliography}{99}
\bibitem{arm85}P. Armitage (1985).  Multistage models of carcinogenesis.  {\it Environmental Health Prespectives} {\bf 63}, 195-201.

\bibitem{armdoll}P. Armitage and R. Doll (1954).  The age distribution of cancer and a multi-stage theory of carcinogenesis. {\it Brit. J. Cancer.} {\bf 8}, 1-12.

\bibitem{armdoll57}P. Armitage and R. Doll (1957).  A two-stage theory of carcinogenesis in relation to the age distribution of human cancer.  {\it Brit. J. Cancer} {\bf 11}, 161-169.

\bibitem{agg}R. Arratia, L. Goldstein, and L. Gordon (1989).  Two moments suffice for Poisson approximations: the Chen-Stein method.  {\it Ann. Probab.} {\bf 17}, 9-25.

\bibitem{beer07}N. Beerenwinkel, T. Antel, D. Dingli, A. Traulsen, K. W. Kinsler, V. E. Velculescu, B. Vogelstein, and M. A. Nowak (2007).  Genetic progression and the waiting time to cancer.  {\it PLoS Comput. Biol.} {\bf 3}, no. 11, 2239-2246.

\bibitem{cal05}P. Calabrese, J. P. Mecklin, H. J. J\"arvinen, L. A. Aaltonen, S. Tavar\'e, and D. Shibata (2005).  Numbers of mutations to different types of colorectal cancer. {\it BMC Cancer} {\bf 5}: 126.

\bibitem{durrsc} R. Durrett (1996).  {\it Stochastic Calculus: A Practical Introduction}.  CRC Press, Boca Raton.

\bibitem{regseq1} R. Durrett and D. Schmidt (2007). Waiting for regulatory sequences to appear.  {\it Ann. Appl. Probab.} {\bf 17}, 1-32.

\bibitem{regseq2} R. Durrett and D. Schmidt (2007). Waiting for two mutations: with applications to regulatory sequence evolution and the limits of Darwinian selection. Preprint, available at http://www.math.cornell.edu/\~{}durrett/recent.html.

\bibitem{DSS}R. Durrett, D. Schmidt, and J. Schweinsberg (2007).  A waiting time problem arising from the study of multi-stage carcinogenesis.  Preprint, available at arXiv:0707:2057.

\bibitem{ek86}S. N. Ethier and T. G. Kurtz (1986).  {\it Markov Processes: Characterization and Convergence.}
John Wiley and Sons, New York.

\bibitem{fh51}J. C. Fisher and J. H. Holloman (1951).  A hypothesis for the origin of cancer foci.  {\it Cancer} {\bf 4}, 916-918.

\bibitem{fn89}D. A. Freedman and W. C. Navidi (1989).  Multistage models for carcinogenesis.  {\it Environmental Health Perspectives} {\bf 81}, 169-188.

\bibitem{harris}T. E. Harris (1963).  {\it The Theory of Branching Processes}.  Springer-Verlag, Berlin.

\bibitem{hknud78}H. W. Hethcote and A. G. Knudson (1978).  Model for the incidence of embryonal cancers: application to retinoblastoma.  {\it Proc. Natl. Acad. Sci. USA} {\bf 75}, 2453-2457.

\bibitem{imkn05}Y. Iwasa, F. Michor, N. L. Komarova, and M. A. Nowak (2005).  Population genetics of tumor suppressor genes.  {\it J. Theor. Biol.} {\bf 233}, 15-23.

\bibitem{imn04}Y. Iwasa, F. Michor, and M. A. Nowak (2004).  Stochastic tunnels in evolutionary dynamics.  {\it Genetics} {\bf 166}, 1571-1579.

\bibitem{ksn03}N. L. Komarova, A. Sengupta, and M. A. Nowak (2003).  Mutation-selection networks of cancer initiation: tumor suppressor genes and chromosomal instability.  {\it J. Theor. Biol.} {\bf 223}, 433-450.

\bibitem{knud71} A. G. Knudson (1971). Mutation and cancer: statistical study of retinoblastoma.
{\it Proc. Natl. Acad. Sci. USA} {\bf 68}, 820--823.

\bibitem{knud01}A. G. Knudson (2001).  Two genetic hits (more or less) to cancer.  {\it Nat. Rev. Cancer} {\bf 1}, 157-162.

\bibitem{kolm38}A. N. Kolmorogov (1938).  On the solution of a problem in biology.  {\it Izv. NII Mat. Mekh. Tomsk. Univ.} {\bf 2}, 7-12.

\bibitem{lm02} E. G. Luebeck and S. H. Moolgavkar (2002).  Multistage carcinogenesis and the incidence of colorectal cancer.  {\it Proc. Natl. Acad. Sci. USA} {\bf 99}, 15095--15100.

\bibitem{mdv88}S. H. Moolgavkar, A. Dewanji, and D. J. Venzon (1988).  A stochastic two-stage model for cancer risk assessment.  I.  The hazard function and the probability of tumor.  {\it Risk Analysis} {\bf 8}, 383-392.

\bibitem{mknud81}S. H. Moolgavkar and A. G. Knudson (1981).  Mutation and cancer: a model for human carcinogenesis.  {\it J. Natl. Cancer Inst.} {\bf 66}, 1037-1052.

\bibitem{ml90}S. H. Moolgavkar and G. Luebeck (1990).  Two-event model for carcinogenesis: biological, mathematical, and statistical considerations.  {\it Risk Analysis} {\bf 10}, 323-341.

\bibitem{ml92}S. H. Moolgavkar and E. G. Luebeck (1992).  Multistage carcinogenesis: population-based model for colon cancer.  {\it J. Natl. Cancer Inst.} {\bf 18}, 610-618.

\bibitem{moran} P. A. P. Moran (1958).  Random processes in genetics.  {\it Proc. Cambridge Philos. Soc.} {\bf 54}, 60-71.

\bibitem{muller} H. J. Muller (1951).  Radiation damage to the genetic material.  In {\it Science in Progress, Seventh Series} (G. A. Baitsell, ed.), Yale University Press, pp. 93-165.

\bibitem{n53} C. O. Nordling (1953).  A new theory on cancer-inducing mechanism.  {\it Brit. J. Cancer} {\bf 7}, 68-72.

\bibitem{nowak} M. A. Nowak (2006).  {\it Evolutionary Dynamics: Exploring the Equations of Life}.  Harvard University Press, Cambridge.

\bibitem{k90}T. Okamoto (1990).  Multi-stop carcinogenesis model for adult T-cell leukemia.  {\it Rinsho Ketsueki} {\bf 31}, 569-571.

\bibitem{sj06}T. Sj\"oblom et. al. (2006).  The consensus coding sequences of human breast and colorectal cancers.  {\it Science} {\bf 314}, 268-274.

\bibitem{wk}D. Wodarz and N. L. Komarova (2005).  {\it Computational Biology of Cancer: Lecture Notes and Mathematical Modeling}.  World Scientific, New Jersey.
\end{thebibliography}
\end{document}